\documentclass{book}

\usepackage{ulem}
\usepackage{url}
\usepackage{epsfig}
\usepackage{graphicx}
\usepackage{makeidx}
\usepackage{multicol}
\usepackage{amssymb}
\usepackage{amsmath}

\usepackage{todonotes}
\graphicspath{ {./Figures/} }
\definecolor{giocolor}{RGB}{0, 150, 100}

\usepackage[T1]{fontenc}
\usepackage[utf8]{inputenc} 			
\usepackage[english]{babel}
\usepackage[percent]{overpic}
\usepackage{comment}
\newtheorem{project}{Project}
\newtheorem{example}{Example}
\newtheorem{exercises}{Exercises}

\makeindex

\def\ttau{\tilde{\tau}}

\newcommand{\pad}[2]{\frac{\partial #1}{\partial #2}}
\newcommand{\tpar}{\tilde{\partial}}

\newcommand{\veps}{\varepsilon}
\definecolor{darkgreen}{rgb}{0.2, 0.5, 0.0}

\usepackage[pageref]{backref}
\renewcommand*{\backrefalt}[4]{%
\ifcase #1 %
(Not cited)%
\or
(Cited on p.~#2)%
\else
(Cited on pp.~#2)%
\fi
}

\usepackage{times}
\usepackage[T1]{fontenc}

\title{Multiscale Modeling with Differential Equations}
\author{Clarissa Astuto, Giovanni Russo}
\date{\today}

\begin{document}

\maketitle
\frontmatter
\tableofcontents

\listoffigures

\mainmatter

\chapter[Multiscale Modeling with Differential Equations]%
{Multiscale Modeling with Differential Equations}
\index{Multiscale!Modeling}

\textbf{Prerequisites}
\begin{itemize}
    \item Differential Equations (DE)     
\end{itemize}
\textbf{Learning Objectives}
\begin{itemize}
\item Learning about mathematical modelling
\item Describe mathematical modeling
\item Show how mathematical techniques can \item Reduce the complexity of a system
\item Define multiple scales
\item Show how the boundary conditions are relevant in a PDE problem
\item Analyze and compare 3D, 2D, and 1D models

\end{itemize}
\textbf{Learning Outcomes} 

Students will learn:
\begin{itemize} 
\item How to solve first order ODE, given the initial condition
\item How to use complex exponentials
\item About techniques based on asymptotic expansions
\item How to use references to build up new knowledge
\item How to practise new knowledge solving the proposed exercises and problems
\end{itemize}

\section*{Contextual Introduction - Summary}
Many physical systems are governed by ordinary or partial differential equations (see, for example, Chapter ''Differential equations'', ''System of Differential Equations''). Typically the solution of such systems are functions of time or of a single space variable (in the case of ODE’s), or they depend on multidimensional space coordinates or on space and time (in the case of PDE’s).

In some cases, the solutions may depend on several time or space scales. An example govered by ODE’s is the damped harmonic oscillator, in the two extreme cases of very small or very large damping (see Section 1.3), the cardiovascular system, where the thickness of the arteries and veins varies from centimeters to microns, shallow water equations, which are valid when water depth is small compared to typical wavelength of surface waves, and sorption kinetics, in which the range of interaction of a surfactant with an air bubble is much smaller than the size of the bubble itself. In all such cases a detailed simulation of the models which resolves all space or time scales is often inefficient or intractable, and usually even unnecessary to provide a reasonable description of the behavior of the system.

In the Chapter ''Multiscale modeling with differential equations'' we present examples of systems described by ODE’s and PDE’s which are intrinsically multiscale, and illustrate how suitable modeling provide an effective way to capture the essential behavior of the solutions of such systems without resolving the small scales. 

\section{Introduction and Examples}
This chapter focuses on space and time multiscale challenges that physical problems face in reality. 

{There are many cases in nature in which problems are characterized by the appearance of structures, whose elements have sizes which span several orders of magnitude. }

\begin{figure}
\centerline{\includegraphics[height=2in]{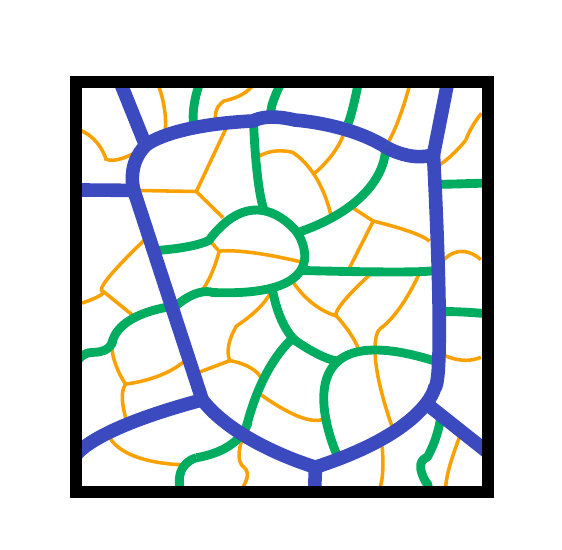}}
\caption{\textit{This figure shows that different sizes of venations are present in the same leaf.}}
\label{fig_leaf}
\end{figure}

{One example that shows up in nature concerns the formation of venation networks in leaves.} The structure of the branches creates a very complicated reticulum, with different ramifications and several different lengths. The lengths of the leaf network can vary from $10^{-5}$ $m$, if we consider veins diameters, to $10^{0}$ $m$, that is the length of some large leaf (see Fig.~\ref{fig_leaf}), {and correspondingly the density of veins dramatically increases as the diameter gets smaller and smaller, thus making it impractical to adopt the same model to describe the flow through large and small veins. See for example \cite{astuto2022comparison,astuto2023asymmetry,astuto2023finite}}. 

When several spatial scales appear in a physical system, for which it is not possible or is impractical to resolve all the spatial scales, it is common to refer to such cases as to \underline{space and/or time}  \underline{multiscale problems} \cite{blonder2020linking}. 

Another well-known example is given by pandemic behaviour. It is a phenomenon that has, not only personal effects, but also, global consequences, such as economic, social and political. To design a good model to describe a pandemic, one has to take multiscale levels into consideration, from individuals to worldwide population, from social distances to social policies, from hand hygiene to the invention of a vaccine that can be distributed to several countries.

In this chapter, the main focus is on four cases, in which different space and time multiple scales are involved. The first application is an example of time multiscale in a damped oscillator. The second topic is  the capture rate of substances freely diffusing around a trap; then we move on with the human cardiovascular system, and, at the end, we describe how to reduce the complexity of hydro flows in channels  . 
\begin{figure}
\centerline{\includegraphics[height=2in]{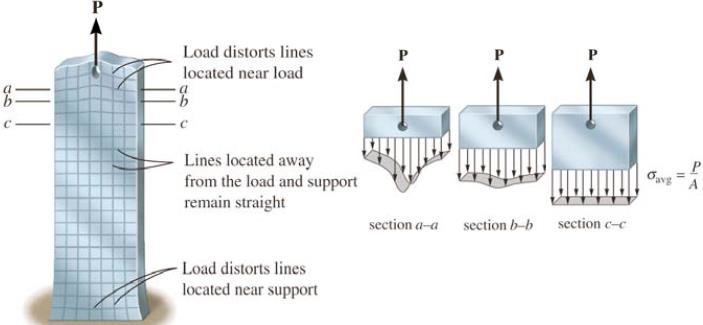}}
\caption{Description of Saint-Venant problem.}
\label{fig_saint_venant}
\end{figure}
Solution of realistic models are often a challenge even for the most advanced computer systems. Very often the complexity comes from the multiscale nature of the processes and the variety of the parameters involved.

Several mathematical models, describing a large range of physical systems, are described by a set of ordinary differential equations (ODEs) or partial differential equations (PDEs), and then solved by suitable numerical methods.

{There is a vast literature on mathematical models described by differential equations, and an equally vast literature on the various techniques for their numerical treatment \cite{credali2022,ASTUTO2023111880,coco2020}.}

Computational power of modern computers is not always sufficient to provide accurate solutions, especially, in three space dimensions, because of the huge number of unknowns that appears in the problem when multiple scales are involved. Therefore, it is necessary to construct suitable multiscale mathematical models of the physical system, that take advantage of the specific nature of the problem, thus allowing a drastic reduction of the number of degrees of freedom (DOFs), but still capturing the essential features of the phenomena that one is interested in.  

{Finding the right boundary conditions can be very helpful in reducing the DOF of a system, when the geometry allows it.} 

In 1855 Adh\'emar Barr\'e de Saint-Venant, a French mathematician and engineer, who gave fundamental contribution in the theory of linear elasticity, introduced the following strategy, known as {\it de Saint Venant's principle} 
to reduce the complexity of a problem in linear elasticity \cite{saint1855memories} with free sides. He ingeniously sidestepped the major difficulty of his problem by replacing  it with a simplified one: consider a cylinder composed of homogeneous isotropic elastic material, with free sides, and in the absence of gravity. The cross sections are arbitrary and, to satisfy static equilibrium, the needed loads are applied, which cause the cylinder to deform, ending up with the \textit{problem of de Saint Venant} in which the goal is to calculate the stress and the displacement throughout the cylinder. According to de Saint Venant principle, 
 the boundary conditions are relaxed, i.e., they are replaced by simpler conditions at the end side of the body, which share with the original one the same total force and torque. In this way, three out of the six components of stress are set to zero, thus drastically simplifying the problem \cite{synge1945problem} (see Fig.~\ref{fig_saint_venant}). When considering torsion or flexion of the cylinder, the problem turns into a boundary condition issue, and {indicates} how to choose the proper conditions to impose (see, for instance, \cite{synge1945problem}). 

This chapter will first show an example within the context of ODE's (see  Section~\ref{sec:oscillator}). We show a simple case in which multiple scales are involved in time, and the parameters of the system vary in a very wide range. We describe two different time scales, $t$, related to the time evolution, and $\tau$, related to the undamped oscillator.   
By examining these two different time scales, we can gain a deeper understanding of how the undamped oscillator behaves and what its characteristics are. Moreover, comprehending how these time scales interact with each other enables us to analyze and make predictions about how the oscillator will behave under various conditions. Some of the applications that could be studied are fascinating phenomena, such as resonance, frequency locking, and other intriguing dynamical effects.

\begin{figure}
\centering
\begin{overpic}[abs,width=0.5\textwidth,unit=1mm,scale=.25]{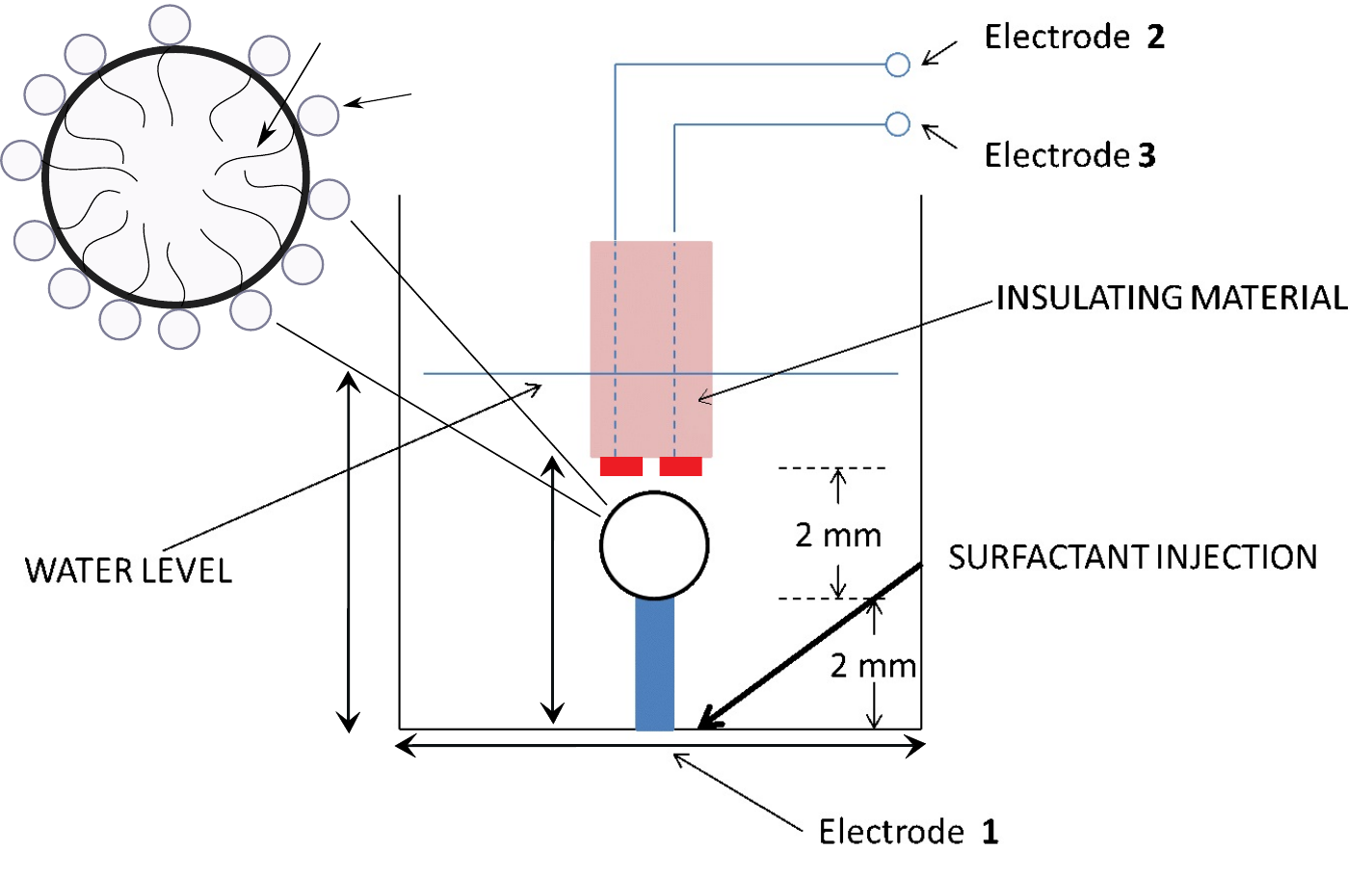}
\put(21,14){{P}}
\put(10,38){\footnotesize hydrophobic tail}
\put(16,36){\footnotesize hydrophilic}
\put(18,33){\footnotesize head}
\end{overpic}
\caption{\textit{Schematic setup of the real apparatus. The central sphere mimics the oscillating gas bubble. The detectors (in red) are located at distance P from the bottom of the vessel.}}
 \label{fig:setup_mm}
\end{figure}

\begin{figure}[tb]
 	\centering
 \includegraphics[width=0.5\textwidth]{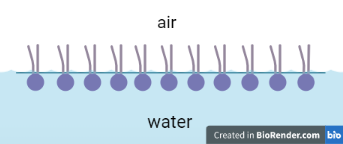}
\caption{\textit{Scheme of surfactant disposition. They are composed by hydrophilic heads, that prefer to stay in the water, and and hydrophobic tails.}}	
	\label{fig:surfactant_2}
\end{figure}

The second application will first focus on \textit{sorption kinetics} application, {in which it is fundamental to couple systems with different dimensions}, such as modeling the capture rate of particles (3D) by moving traps (2D). This is a common research topic among different fields of sciences, such as in chemistry, when studying biomolecular reactions, or in biology, with {complex prey-predator ecological systems.} The last example is given by living cells and their moving
membranes. The nature of this motion is associated to the exchange of thermal and metabolic energy \cite{brochard1975frequency,Raudino20168574,dukhin1995dynamics,nguyen2013behavior} and it has been studied how the oscillations of the cell membrane affect the capture rate of substances freely diffusing around them. To obtain quantitative results, a reproducible and tunable biomimetic experimental {setup} has been considered to simulate the phenomenon. 
The system consists of an oscillating air bubble in water and a diffusive flow of ions, see Fig.~\ref{fig:setup_mm}. When the negative ions meet the surface of the bubble, they are partially adsorbed because of their polar head, but at the same time, they also prefer to settle at the bubble surface because of the hydrophobic repulsion between apolar tails and the water where it is immersed (see Fig.~\ref{fig:surfactant_2}). This sort of substance is called surfactant, and in Section~\ref{sec:sorption} we see it in details.

Modeling diffusion and trapping in presence of a bubble is very challenging because of the multiple scales involved in space. The range of the attractive-repulsive core of a trap is of the order of tens of nanometers ($10^{-9}m$),  a length that is prohibitively small for detailed simulations (see , e.g., \cite{CiCP-31-707}). 
while the size of the domain is a few millimeters ($10^{-3}m$), thus almost 6 different orders of magnitude in space.

To overcome such difficulties, in \cite{multiscale_mod} the authors propose a \emph{multi-scale model} to effectively treat the problem describing the trapping by a suitable boundary condition on the bubble surface. The potential is assumed to have a range of small size $\varepsilon$. An asymptotic expansion in $\varepsilon$ is considered, and the boundary conditions are obtained by retaining the lowest order terms in the expansion. 
{In the reduced multiscale model,}  the ions are {described} by a volumetric concentration which satisfies a simple diffusion equation in the domain of the liquid (3D), while the carriers that are trapped on the surface satisfy an equation for the superficial concentration (2D), that is coupled with the parabolic one. The two concentrations are related by a mass balance equation which acts as a suitable boundary condition for the concentration in the bulk. 
{This model is described in detail in Section~\ref{sec:sorption}.}

In Section~\ref{sec:cardiovascular} we refer to \textit{geometric multiscale} approach for the numerical simulation of blood flow problems, a very hot topic in applied mathematics (see, e.g.,  \cite{quarteroni2016geometric,quarteroni2017cardiovascular}). The geometric multiscale analysis refers to the idea of coupling different dimensions in order to obtain a detailed study of the phenomenon. It plays an important role when different scales or levels interact in a system.
With this approach, different features of the entire cardiovascular system can described, focusing on the right scale that is more suitable for the problem.
\begin{figure}
\centerline{\includegraphics[height=2in]{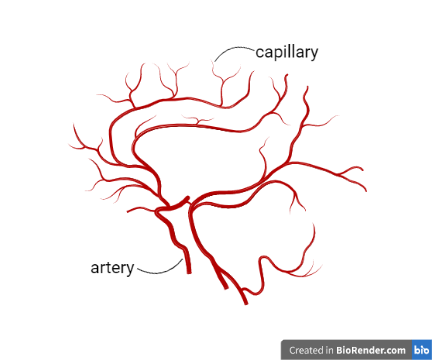}}
\caption{\textit{Representation of large arteries and thin capillaries in the circulatory system.}}
\label{fig_capillaries_veins}
\end{figure}

{Human circulatory system is well approximated by hyperbolic equations. Even if blood can be considered an incompressible fluid, the elasticity of veins and arteries introduces compressibility in the model, with finite speed of propagaton, which is much smaller than the propagation of sound waves in water. The first one who developed a mathematical model for liquid flow in an elastic duct was Leonard Euler, as appears in a paper published after his death \cite{EulerBlood}\footnote{It is conceivable that he wrote the paper  even before his celebrated papers on Euler equations of compressible gas dynamics, datin 1757 \cite{euler1757principes}}. 
In this respect he may be considered the father of hemodynamics (as well as of innumerable other mathematical discoveries!)}

{Let us restrict, for example, to the pulmonary blood circulatory system.} The vessels, have a wide range of different diameters. Pulmonary arterial and vein trees are composed by many orders of branches, with diameters ranging from few $\rm \mu$m for capillaries to about 2-3 cm for the main pulmonary artery \cite{townsley2012structure}.

This is what is called space multiscale challenges. Trying to solve so many different orders of magnitude with the same space discretization, would be very computationally expensive and, in certain cases, not effective, unless all space scales are very well resolved. Mathematical modeling and approximation are needed to investigate such cases.

\begin{figure}
\centerline{\includegraphics[height=2in]{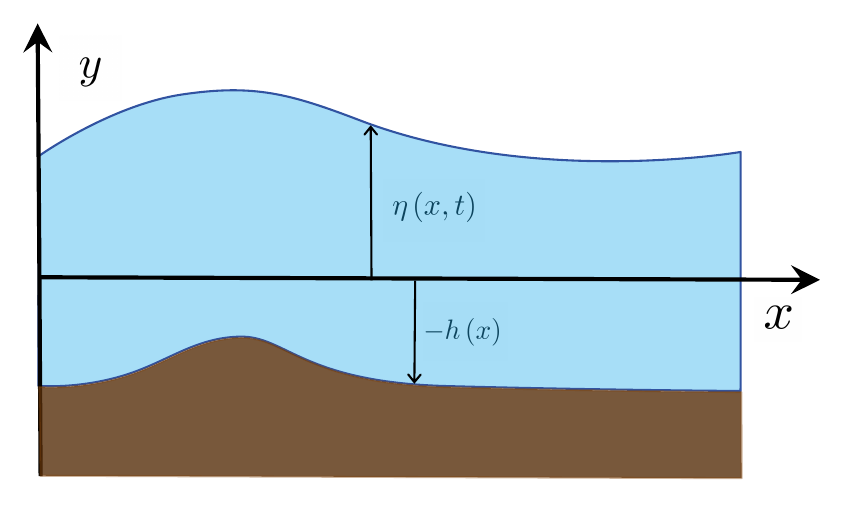}}
\caption{Shallow water equation is an approximation of NS equations, when the horizontal length scale is much greater than the vertical length scale, i.e. $\eta + h \ll $ x-direction of the flow.
}
\label{shallow_water1}
\end{figure}

{Going back to the idea of reducing the number of DOF, a common technique in such cases consists in adopting different spatial models in different parts of the system, and coupling them together. 
In particular, in regions where high accuracy is not required, one can adopt a lower dimensional model.} Practically, the geometry of the problem reveals the region where the study is focused on, describing it with a 3D model. The description of the background is modeled with lower dimensional problems. As a result, the problem is solved coupling 3D - 1D, 3D - 0D or 1D - 0D systems. A three-dimensional (3D) model can be reduced to a one-dimensional (1D) model, or to a zero-dimensional (0D) model, where only the dependence on time survives. 

{Matching domains with different dimensions  can generate inconsistency at the boundaries. For this reason, the boundary conditions play an important role in the system. In particular, they play a key role in determining suitable treatment of the junctions, i.e., those parts of the system which separate domains of different dimensions. A similar problem appears in systems that treat water channels (Section 1.3) and gas tubes. }

{In summary, the simulation of the entire cardiovascular system of human body is a challenging multiscale problem. Adopting a multiscale approach, by  coupling domains of different size using different spatial dimensions, allows  numerical simulations of cardiovascular diseases, thus describing many features of medical interest with progressively more accuracy and precise quantification of the error.}

The last example of space scales reduction is shown in Section~\ref{sec:SW_channels}. The coupling of different dimensions is also widely used to describe flows in artificial channels with applications to environmental problems \cite{bellamoli2018numerical,toro_ADER,toro_WAF}, where the basic model comes from the Shallow Water (SW) equations. These are approximations of incompressible, irrotational Euler equation with free boundary, valid when the depth of the water is small compared to the wave-length. Under this assumption it can be shown that the pressure is basic the hydrostatic pressure, i.e.\ in each point it is proportional to the height of the column of water above that point. Section~\ref{section_sw} in the Appendix  is dedicated to the extended derivation of the SW equations. In the 1D case such equations are also called de Saint-Venant equations (yes, the same de Saint Venant of linear elasticity!) \cite{SV-shallow-water}.

\begin{figure}
\centerline{\includegraphics[width=0.45\textwidth]{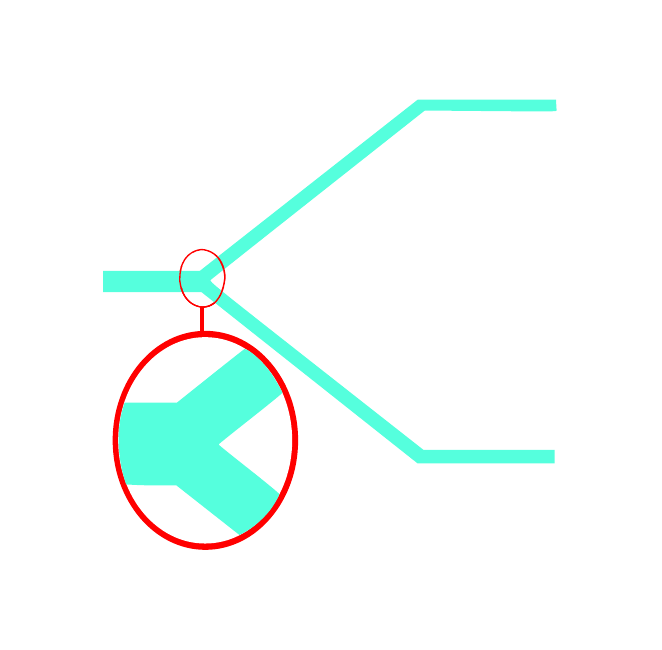}}
	\caption{\textit{Example of a channel network. Regions with two dimensional behaviour are zoomed in.}}
	\label{fig:junction}
\end{figure}

In nature, many phenomena are described by SW equations, such as flows in open channels and rivers. It has been proposed a 2D - 1D model for the junction of shallow-water channel networks, with a localized 2D problem at the junctions and a 1D model for the channels, because it is much more efficient from a computational point of view. The difficulty is in finding the right expression for the conditions at the junction between the channels; see Fig.~\ref{fig:junction}.

\section{Time multiscale in a damped oscillator}
\label{sec:oscillator}
Multiscale problems can appear in time as well, when very different time scales are involved in a phenomenon. 

Let us consider again the example of the damped oscillator, seen in Chapter \textbf{Dimensional Analysis}. The equations governing its motion are given by:
\begin{align}
\label{eq_xpunto}
    \frac{dx}{dt} &= v \\
\label{eq_vpunto}
    m\frac{d v}{dt} &= -kx - \gamma v
\end{align}
where $m,k$ and $\gamma$ are, respectively, the mass of the oscillator, the elastic constant of the spring and the damping coefficient.

Eqs.~(\ref{eq_xpunto}-\ref{eq_vpunto}) are a system of two differential equations of first order. Differentiating the first equation \eqref{eq_xpunto} and making use of the second one \eqref{eq_vpunto}, the dynamics can be described by a differential equation of second order:
\begin{equation}
\label{eq_x2punti}
    \frac{d^2 x}{dt^2} = -\frac{\gamma}{m} \frac{dx}{dt} - \frac{k}{m}x.
\end{equation}
Eq.~\eqref{eq_x2punti} is an example of ordinary differential equation of degree two, with constant coefficients.

If $k=0$, there is no spring, and the equation for $v$ reduces to
\begin{equation}
\label{eq_v}
    \frac{dv}{dt} = -\frac{\gamma}{m}v.
\end{equation}
The solution can be found by observing that the time derivative of the function $v$ is proportional to itself. The functions that are proportional to their own derivatives are the exponentials, therefore it is natural to look for a solution of the form
\begin{equation}
\label{eq_exact_exp}
    v = e^{\lambda t}.
\end{equation}
By inserting Eq.~\eqref{eq_exact_exp} into Eq.~\eqref{eq_v}, one obtains
\begin{equation}
    \lambda e^{\lambda t} = -\frac{\gamma}{m} e^{\lambda t}
\end{equation}
from which we can deduce that 
\begin{equation}
    \lambda = -\frac{\gamma}{m}
\end{equation}
thus, one solution of Eq.~\eqref{eq_v} is 
\begin{equation}
   \displaystyle \bar v(t) = e^{-\frac{\gamma}{m}t}.
\end{equation}

Since the equation for $v$ is linear, if a given function $\bar v(t)$ is a solution, then $C\bar v(t)$ is also a solution, with $C$ is an arbitrary constant. If we impose that such a function satisfies the initial condition $v(0) = v_0$, then this will fix the constant $C$, such that $C = v_0$, and the solution will be
\begin{equation}
    v(t) = v_0 e^{-\frac{\gamma}{m}t}.
\end{equation}

From Eq.~\eqref{eq_xpunto} it is possible to compute the position of the mass by integrating the velocity:
\begin{align}
    x(t) & = x_0 + \int_0^t v(t')dt' \nonumber \\
         & = x_0 + v_0\int_0^te^{-\frac{\gamma}{m}t'} dt' \nonumber \\
         & = x_0 - v_0\frac{m}{\gamma}\left[ e^{-\frac{\gamma}{m}t'}\right]_0^t \nonumber \\
         & = x_0 - v_0\frac{m}{\gamma}\left( e^{-\frac{\gamma}{m}t} -1\right),
\end{align}
where $x_0 = x(0)$ and, if $t\to \infty$, then $x\to x_{\rm f} = x_0 + \frac{v_0\, m}{\gamma}$. 

Notice that, dimensionally, $\displaystyle  \left[ \frac{\gamma}{m}\right] = T^{-1}$, so $\displaystyle \frac{m}{\gamma}$ is the decay time due to the damping, which we denote by $t_d = \frac{m}{\gamma}$.

This result shows that, with no spring a mass would slow down to stop at final position $x_{\rm f}$, starting at $x_0$ with initial velocity $v_0$.

Now let us consider the case in which $\gamma = 0$. In such cases the equation for the position $x(t)$ would be
\begin{equation}
\label{eq_2nd_ord}
    \frac{d^2 x}{d t^2} = -\frac{k}{m}x.
\end{equation}
This equation tells us that the function $x(t)$ is such that it is proportional to its second derivative, but with a different sign. What are the functions whose second derivative are proportional to the opposite of themselves? Those are the trigonometric functions! So we look for solutions 
 of the form
\begin{equation}
    x(t) = \sin(\omega t) \quad {\rm or}\quad  x(t) = \cos(\omega t).
\end{equation}
In the first case we have
\begin{align}
    x'(t) &= \omega \cos(\omega t), \\ 
    x''(t) &= -\omega^2\sin(\omega t) = -\frac{k}{m}x(t) = -\frac{k}{m}\sin(\omega t)
\end{align}
indicating that $\sin(\omega t)$ is a solution if $\omega^2 = {k}/{m}$. Likewise, in the second case we have that $x(t) = \cos(\omega t)$ is a solution, still providing $\omega^2 = {k}/{m}$. 

Since Eq.~\eqref{eq_2nd_ord} is linear, if $\sin(\omega t)$ and  $\cos(\omega t)$ are solutions, then any linear combination is a solution. The general solution is therefore
\begin{align}
\label{eq_general_solu}
    x(t) &= a\cos(\omega t) + b\sin(\omega t) \\
    v(t) &= -a\,\omega \sin(\omega t) + b\,\omega \cos(\omega t).
\end{align}
The coefficients $a$ and $b$ can be determined by imposing the initial conditions 
\begin{equation}
\label{ic_trigonom}
x(0) = x_0, \quad v(0) = v_0.
\end{equation}

Observe that the expression of the two solution is not so different if we make use of complex exponentials. Indeed, since we have Euler formula
\begin{equation}
    e^{i\omega z} = \cos(\omega t) + i \sin(\omega t)
\end{equation}
and its inverse:
\begin{equation}
    \cos(\omega t) = \frac{e^{i\omega t} + e^{-i\omega t}}{2}, \quad
    \sin(\omega t) = \frac{e^{i\omega t} - e^{-i\omega t}}{2i}
\end{equation}
where $i$ denotes the imaginary units, i.e., $i^2 = -1$, then the solution of the undamped harmonic oscillator can be written in terms of complex exponentials:
\begin{equation}
\label{eq_gener_solu_complex}
    x(t) = C_+ e^{i\omega t} + C_- e^{-i\omega t}, \qquad \omega = \sqrt{\frac{k}{m}} = \frac{1}{t_0}.
\end{equation}
The time constant $t_0 = \omega ^{-1}$ is a related to a full period $T$ of the undamped oscillator by $T = 2\pi t_0$, because after one period the argument has changed by $\omega T=2\pi$.

Here we focus on the general case, i.e., $k\neq 0, \gamma \neq 0$. Eq.~\eqref{eq_x2punti} for the position has the following structure
\begin{equation}
    \label{eq_x2punti_general}
    x'' + \widetilde \gamma x' + \omega^2 x = 0, \qquad \widetilde \gamma = \frac{\gamma }{m} = \frac{1}{t_d}.
\end{equation}
where for short we denote by $~'$ the time derivative. It is therefore natural to look for solutions of the form  $x(t) = e^{\lambda t}$, and inserting this expression in Eq.~\eqref{eq_x2punti_general} we obtain that $e^{\lambda t}$ is indeed a solution provided 
\begin{equation}
\label{eq_lambda_2nd}
    \lambda^2 + \widetilde \gamma \lambda + \omega^2 = 0.
\end{equation}
This is an algebraic equation of degree 2. The discriminant is given by $\Delta = \widetilde \gamma^2 - 4\omega^2$.  If $\Delta > 0 $ there will be two real negative roots, while if $\Delta < 0$ there will be two complex conjugate roots.

If $\Delta < 0$, i.e.\ $\widetilde \gamma < 2\omega$ or, equivalently, if $t_d>2t_0$, the solution shows damped oscillations, if $\Delta > 0$, i.e.\ if $t_d<t_0$, the solution is damped with no oscillations, while if $\Delta = 0$ we have the so called critical damping\footnote{For a given fixed oscillation frequency $\omega = \sqrt{k/m}$, the value of damping $\widetilde \gamma^*$, that gives $\Delta = 0$, is the one which has the fastest damping.}.

\subsection{Dimensional analysis}
The Eqs.~(\ref{eq_xpunto}-\ref{eq_vpunto}), governing the damped oscillator, apparently depend on three parameters, $m,k$ and $\gamma$. However, we observe that after dividing by the mass $m$, there are only two parameters the equation depends on. If we use non-dimensional time to describe the evolution of the system,  we can show that there is only one parameter that determines the quantitative behaviour of the solution. 

Indeed, in an undamped oscillator, the quantity $\omega t$ can identify a non dimensional time $\tau$. In term of this parameter the equations become:
\begin{align}
\label{eq_xpunto_tau}
    \omega \frac{dx}{d\tau} &= v \\
    \omega \frac{d v}{d\tau} & = -\omega^2x - \frac{\gamma}{k}v. \label{eq_vpunto_tau}
\end{align}
Differentiating \eqref{eq_xpunto_tau}, and replacing the expression of $dv/d\tau$ in the second equation, after dividing by $\omega^2$ we obtain
the following second order differential equation for $x$:
\begin{equation}
    \frac{d^2 x}{d\tau^2} + f\frac{d x}{d\tau} + x = 0
\end{equation}
and the discriminant is $\Delta = f^2 - 4$, so the critical damping occurs if $ f = 2$. If $f < 2$, i.e., $t_0 < 2t_d$, there are damped oscillations, and if $f > 2$ there are no oscillations. If $f \to \infty$, the equation for $v$ suggests that $v \to 0$, and $x$ is constant, so the limit of infinite damping is not very interesting.

If we want to observe a non trivial behaviour for $f \to \infty$, we have to act differently. 

Let us start again from the original system (\ref{eq_xpunto}-\ref{eq_vpunto}), and let us divide the second equation by $k$, obtaining 
\begin{equation}
    \frac{m}{k}\frac{dv}{dt} = -x -\frac{\gamma}{k}v.
\end{equation}
We can write this equation as 
\begin{equation}
    t_0^2\frac{dv}{dt} = -x-t_\gamma v
\end{equation}
where $\displaystyle t_0 = \sqrt{\frac{m}{k}}$ and $\displaystyle t_\gamma = \frac{\gamma}{k}$. Notice that $t_d\,t_\gamma = t_0^2$.

We now define a {\em new\/} non-dimensional time $\ttau$ as $\displaystyle \ttau = \frac{t}{t_\gamma}$, with this definition, Eq.~\eqref{eq_xpunto_tau} becomes:
\begin{equation}
    \frac{t_0^2}{t_\gamma^2}\frac{d^2x}{d\ttau^2} = -x - \frac{dx}{d\ttau}.
\end{equation}
As $\displaystyle \frac{t_0}{t_\gamma} \to 0$ we obtain the equation 
\begin{equation}
    \frac{dx}{d\ttau} = -x
\end{equation}
whose solution is $x = x_0e^{-\ttau} = x_0e^{-t/t_\gamma}$.

Note that, the limit is different form the previous one. Indeed, the original second order system (\ref{eq_xpunto}-\ref{eq_vpunto}), which we can rewrite in the form 
\begin{align}
    \frac{dx}{dt} & = v \\
    t_0^2 \frac{dv}{dt} & = -x -t_\gamma v
\end{align}
will relax to the system
\begin{align}
    \frac{dx}{dt} &= v \\
    x + t_\gamma v & = 0,
\end{align}
as $\ttau_0 \to 0$, which is a first order system of ODE's.

The solution of such system requires fewer initial conditions: the initial condition of the position, $x_0$, is sufficient, and the velocity is related to the position by the algebraic equation $x + t_\gamma v = 0$.

\subsection{Singular perturbation problems and differential-algebraic equations}
This example is a particular case of a large class of systems in which there is a small parameter that introduces a small scale. Such systems take the form:
\begin{align}
    \frac{dx}{dt} &= f(t,x,y) \\
    \varepsilon \frac{dy}{dt} & = g(t,x,y)
\end{align}
where $x(t) \in \mathbb{R}^{n_1}$, $y(t) \in \mathbb{R}^{n_2}$.

As $\varepsilon \to 0$, the system degenerates into a system of $n_1$ differential equations, and $n_2$ algebraic equations of the form
\begin{align}
    \frac{dx}{dt} &= f(t,x,y) \\
      0 & = g(t,x,y).
\end{align}
This kind of systems are called differential-algebraic systems (DAE) (see, e.g. \cite{ascher1998computer}).

\begin{figure}[!ht]
	\centering
\begin{overpic}[abs,width=0.5\textwidth,unit=1mm,scale=.25]{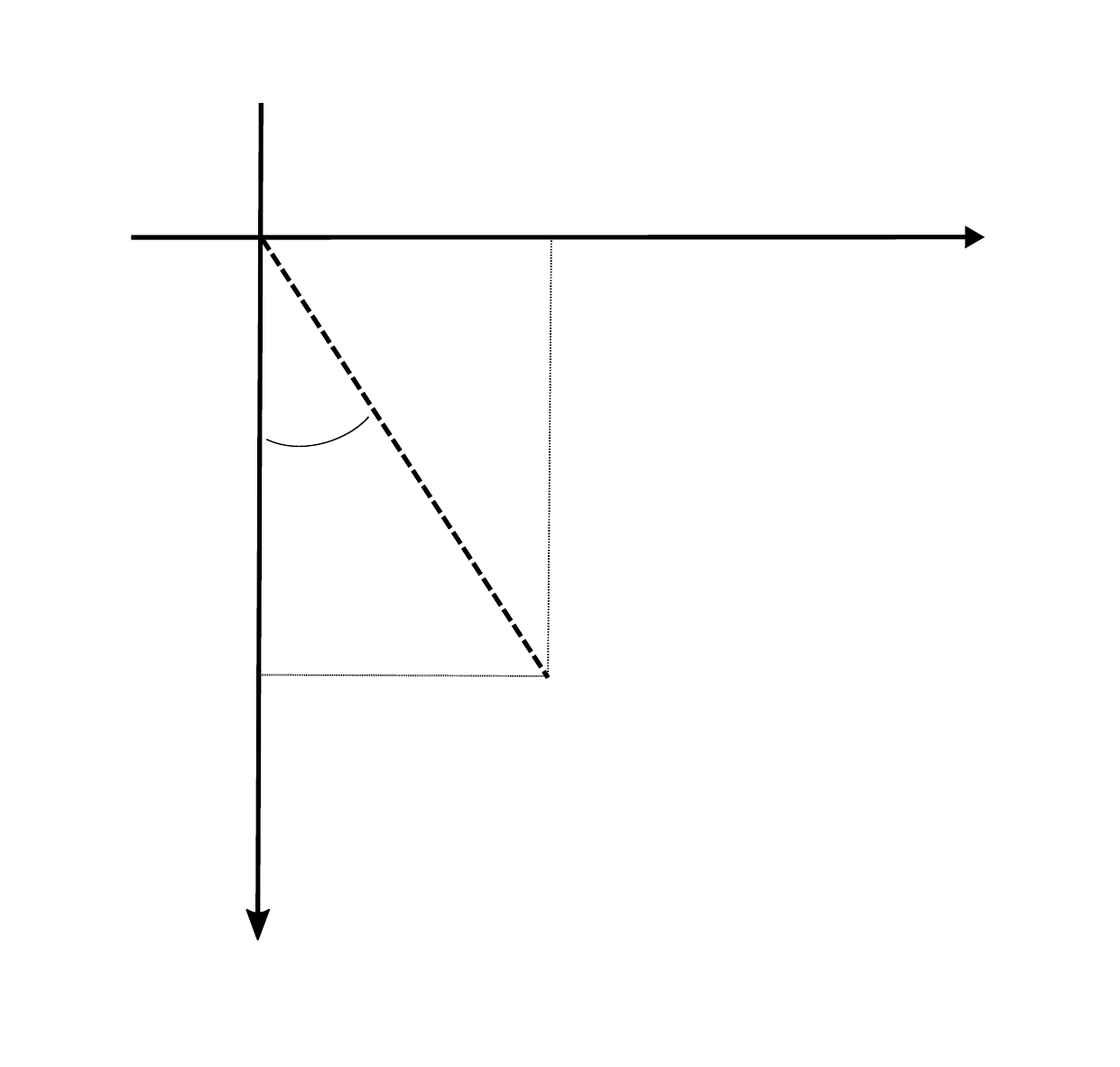}
\put(17,30){$\theta$}
\put(31,20){$m$}
\put(9,21){$x_2$}
\put(30,46){$x_1$}
\end{overpic}
	\caption{\textit{Schematic representation of system~(\ref{eq_m_x12punti}-\ref{eq_m_x22punti}).}}
	\label{fig_axis_x1_x2}
\end{figure}

An example of such system is given by the pendulum obtained by a spring, with elastic constant $k$ and rest length $L$, fixed on one end, and with a mass $m$ on the other end (see Fig.\ref{fig_axis_x1_x2}). The equations of motion of the coordinates of the system are given by Newtons second law:
\begin{align}
\label{eq_m_x12punti}
    m \ddot{x}_1 &= -k (x_1-L\sin\theta) \\
\label{eq_m_x22punti}
    m \ddot{x}_2 &= -k (x_2-L\cos\theta)+mg 
\end{align}
Here
\begin{equation}
    x_1 = L\sin\theta, \quad x_2 = L\cos\theta \label{eq:x1_x2}
\end{equation}
from which we can deduce 
\[
    \sin \theta = \frac{x_1}{\sqrt{x_1^2 + x_2^2}},\, \cos \theta = \frac{x_2}{\sqrt{x_1^2 + x_2^2}},
\]
USing these relations, the equations of motion can be written as:
\begin{align}
    m \ddot{x}_1 &= -k x_1\left(1-\frac{L}{\sqrt{x_1^2 + x_2^2}}\right) \\
    m \ddot{x}_2 &= -k x_2\left(1-\frac{L}{\sqrt{x_1^2 + x_2^2}}\right)+mg 
\end{align}
As $k\to \infty,$ then $\displaystyle k\left(1-\frac{L}{\sqrt{x_1^2 + x_2^2}}\right)$ becomes undetermined, since $k\to \infty$ and \\ $\displaystyle \left(1-\frac{L}{\sqrt{x_1^2 + x_2^2}}\right)\to 0$. Let $\lambda$ denote such limit. Then the equations of motion become
\begin{equation}
    m \ddot{x}_1 = -\lambda x_1, \quad
    m\ddot x_2 = -\lambda x_2 + mg \label{eq:pendulum2}
\end{equation}
and $\lambda$ is determined by imposing that $x_1^2 + x_2^2 = L^2$.
The solution of this system provides not only the dynamics of the mass, but also the tension $\lambda$ of the rigid spring. 
It is possible, with suitable manipulations, the classical equations of motion of the pendulum
\begin{equation}
    L\frac{d^2\theta}{dt^2}+g \sin\theta = 0 \label{eq:pendulum1}
\end{equation}
starting from equations \eqref{eq:pendulum2} and \eqref{eq:x1_x2} (see Exercise \ref{ex:pendulum}).

\begin{figure}[tb]
	\begin{minipage}
		{.3\textwidth}
		\centering
		\includegraphics[width=1.1\textwidth]{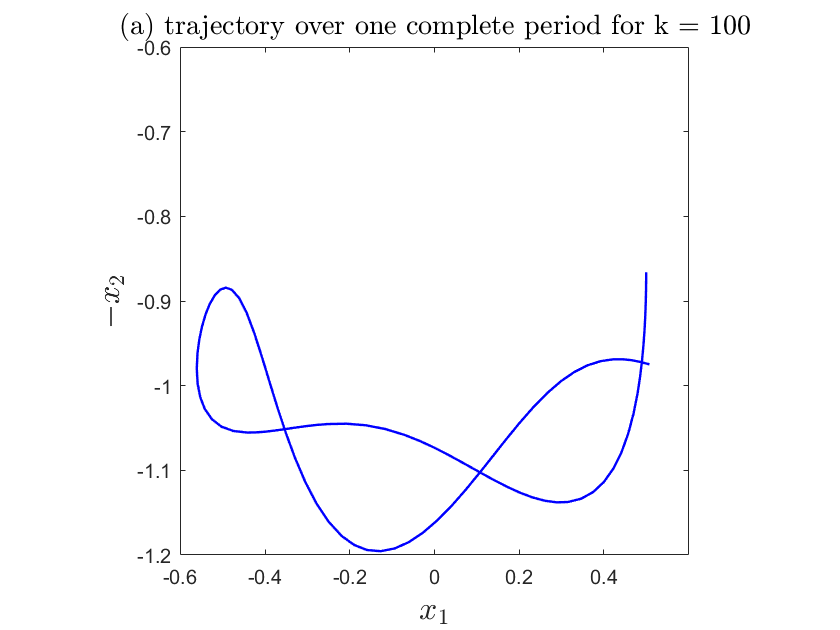}
	\end{minipage}
	\begin{minipage}
		{.3\textwidth}
		\centering
		\includegraphics[width=1.1\textwidth]{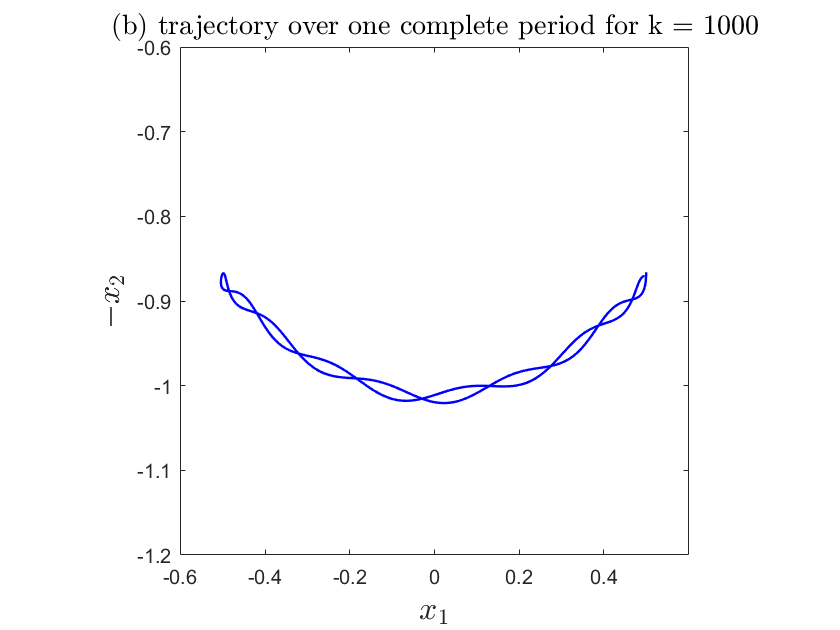}
	\end{minipage}
	\begin{minipage}
		{.3\textwidth}
		\includegraphics[width=1.1\textwidth]{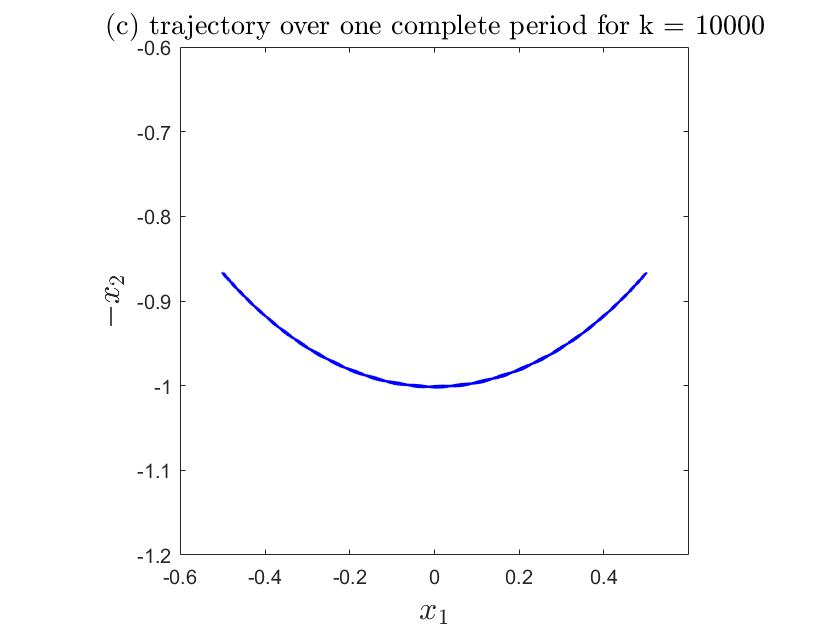}
	\end{minipage}
	\caption{\textit{In these tree plots we see  different oscillations for different elastic constants, $k$, namely $k=10^2, 10^3, 10^3$ (N/m). When $k\to \infty$, the oscillator behaves like a classical pendulum.}}
\label{fig_pendulum}
\end{figure}

Fig.~\ref{fig_pendulum} shows the trajectory of the mass of 1kg, oscillating in a pendulum with rest length $L = $ 1m, starting from an angle $\theta_0=30^\circ$, with an initial length $L$, for various values of the stiffness parameter. 
As the stiffness $k$ of the spring increases, the trajectories resemble more and more those of a classical pendulum.

\section{Coupling conditions for sorption kinetics at the oscillating cell membrane}
\label{sec:sorption}

The term surfactant comes from \textbf{surf}ace \textbf{act}ive \textbf{a}ge\textbf{nt}. They are composed by anions and cations, that have two different configurations. The negative ions are amphiphilic molecules, i.e., they both have hydrophilic and hydrophobic parts, and, for this reason, they are absorbed at the air-water interface. The hydrophobic part prefers to stay in the air, while the hydrophilic part stays in water. This will cause a decrease in the surface (or interfacial) tensions.  



In this section we see how to treat different orders of magnitude in space, in the same mathematical model. How we were saying in the Introduction, an expansion in $\varepsilon$ (the range of the attractive-repulsive potential) will be considered, in order to obtain a new boundary condition for the reduced multiscale model. 

Introducing the local concentration of ions $c= c(\vec{x},t)$, its time evolution in a fluid  is governed by following the conservation law
\begin{align} 	\label{equation_flux}
	\displaystyle \frac{\partial c\left(\vec{x},t\right)}{\partial t} &=-\nabla\cdot J(\vec{x},t) \qquad {\rm in }\, \Omega \\
		\displaystyle J &=\ -D\left(\nabla c +\ \frac{1}{k_BT}c\nabla V\right) \qquad {\rm in }\, \Omega   		\label{equation_definition_flux}
\end{align}
where $k_B$ is the Boltzmann's constant, $T$ is the absolute temperature (assumed to be constant) and $V(x)$ is a suitable potential function that simulates the \textit{attractive-repulsive} behavior of the bubble with the particles. 

A part from the attraction, the potential models also the impermeability of the bubble, thus, at very short distances, it is designed to repulse the particles from the surface of the trap, and, for this reason, its derivative becomes negative in the proximity of the bubble radius.


The multiscale model has been derived and tested in various dimensions $d$. 

\subsection{1D Multiscale model}
\label{section_1Dmodel}
\begin{figure}[!ht]
	\centering
	\begin{minipage}
		{.49\textwidth}
		\centering
		\includegraphics[width=\textwidth]{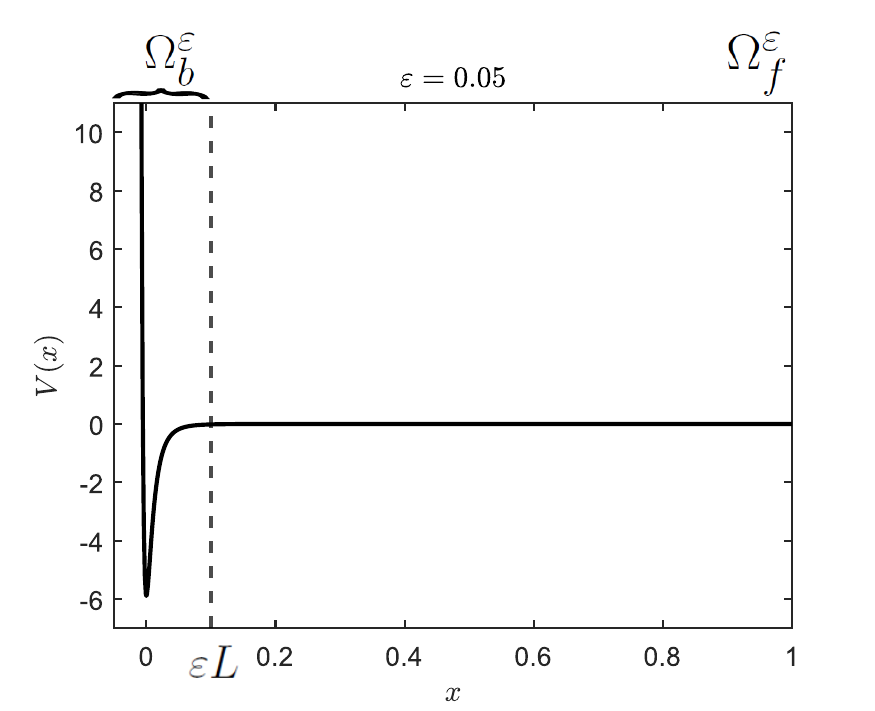}
	\end{minipage}
	\begin{minipage}
		{.49\textwidth}
		\includegraphics[width=\textwidth]{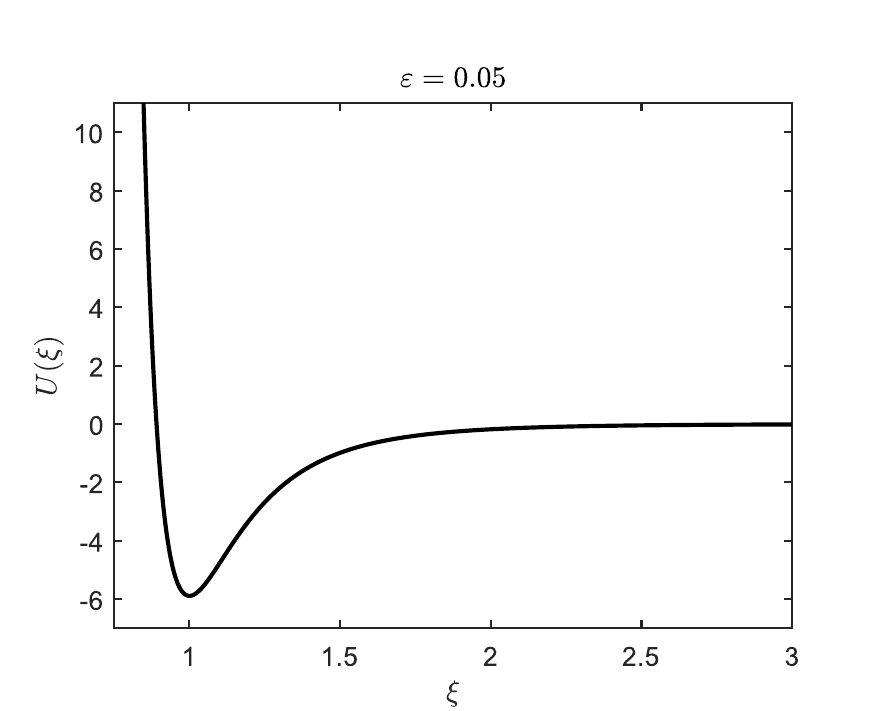}
	\end{minipage}	\caption{\textit{Representation of the external potentials $V(x)$ (in Eq.~\eqref{expr_V_LJ}), on the left, and $U(\xi)$ (in Eq.~\eqref{expr_U_LJ}), on the right, for $\varepsilon = 0.05$. $U(\xi)$ is obtained from a change of variable of $V(x)$, with $\xi \in [0,L+1]$. }}
	\label{figure_potential_V_1D}
\end{figure}

This subsection is devoted to the derivation of the multiscale model in 1D space dimension. 

The drift-diffusion Eqs.~(\ref{equation_flux}-\ref{equation_definition_flux}) in 1D read
\begin{eqnarray}
	\label{pde1d}
	\displaystyle \frac{\partial   c }{\partial t} + \frac{\partial  J}{ \partial x} = 0 \qquad {\rm in} \, \Omega
	\end{eqnarray}
where
\begin{equation}
\label{flux1d}
	\displaystyle J = - D\left( \frac{\partial  c  }{ \partial x} + \frac{1}{k_BT}  c  V'\right),
\end{equation}
and {  the full domain is $\Omega = [-\varepsilon,1]$, with $\Omega^\veps_b = [-\veps,L\veps]$ the bubble domain, inside which the ions feel the effect of the potential, and} the fluid domain is $\Omega^\varepsilon_f = [L\varepsilon,1],$ which is not affected by the bubble, see Fig.~\ref{figure_potential_V_1D}. 
Therefore, the potential has effect only in $\Omega^\varepsilon_b$, while it is constant for $x \in \Omega^\varepsilon_f$, say $V(x)=0$. 

A suitable potential is the Lennard-Jones (LJ) type, which describes attraction at long distances and repulsion at short distances, due to Van der Waals and Pauli terms, respectively. A typical shape of the potential $V(x)$ is shown in Fig.~\ref{figure_potential_V_1D} 
and it takes the form
\begin{eqnarray}
\label{expr_V_LJ}
	V(x) = E\left( \left(\frac{x+\varepsilon}{\varepsilon}\right)^{-12} - 2\left(\frac{x+\varepsilon}{\varepsilon}\right)^{-6} \right)
\end{eqnarray}
where $\varepsilon$ denotes the range of the potential and $E$ represents the depth of the well (see Fig.~\ref{figure_potential_V_1D}, left panel). It is convenient to adopt a non dimensional form of the potential, expressed as a function of the rescaled variable $\xi = 1 + x/\varepsilon \in [0,L+1]$:
\begin{equation}
\label{expr_U_LJ}
    U(\xi) = \phi \left( \xi^{-12} - 2\xi^{-6} \right)
\end{equation}
where { $\displaystyle \phi = {E}/{k_BT}$ represents the ratio between the depth, $E$, of the potential well of the bubble and $k_B T$. 
$L$ is a scaled distance beyond which the potential is negligible. 
}


The assumption made is that the trap is located at $x=0$ and its ''thickness'' is of order of $\varepsilon$.
The original model requires to solve Eq.~\eqref{pde1d} in $\Omega^\varepsilon=\Omega^\varepsilon_b \cup \Omega^\varepsilon_f$ with initial condition $ c  =  c_0(x)$ when $t=0$ and with boundary conditions
\begin{equation}\label{bc1d}
	J=0 \quad {\rm at} \> x=-\varepsilon, {
	\rm and} \> x = 1
\end{equation}
to simulate wall effects. 

It is clear that the solution $ c $ and the flux $J$ depend on $\varepsilon$ as well, and then the original problem can be written as
\begin{eqnarray}
\label{eq_flux_eps}
\displaystyle \frac{\partial  c _\varepsilon}{\partial t} + \frac{\partial J_\varepsilon}{\partial x} &=&0\quad \rm{in } \, \Omega^\varepsilon\\ \label{eq_flux_eps_2}
\displaystyle J_{\varepsilon} &=& -D\left(\frac{\partial  c _\varepsilon}{\partial x} + \frac{1}{k_B T} \,  c _\varepsilon  \, V'_\varepsilon \right).
\end{eqnarray}
It is important to notice that the space derivative can be expressed ad $\partial / \partial x = \partial /\partial \xi \cdot \partial \xi / \partial x = \partial /\partial \xi \cdot 1/ \varepsilon$, and, with this notation, Eq.~\eqref{eq_flux_eps_2} becomes
\begin{equation} \label{newflux}
J_\varepsilon = -D \frac{1}{\varepsilon}\left(\frac{\partial  c _\varepsilon}{\partial \xi}+ c _\varepsilon \, U' \right).  
\end{equation}
{Rewriting \ref{eq_flux_eps} in the rescaled variable and using Eq.~\eqref{newflux}, it becomes
\begin{equation}
    \label{c_eps2}
    \frac{\partial  c_\varepsilon}{\partial t} = -\frac{1}{\veps}\pad{J_\veps}{\xi} = \frac{D}{\varepsilon^2}\frac{\partial}
    {\partial \xi} \left(\pad{c_\veps}{\xi} + c_\veps U'(\xi)\right).
\end{equation}
}

The range of the scaled variable $\xi$ does not depend on $\varepsilon$, then the following expansion in $\Omega^\varepsilon_b$ is valid for the concentration $ c _\varepsilon(\xi,t)$  :
\begin{equation}\label{exprho}
{c_\varepsilon(\xi,t) =  c ^{(0)}(\xi,t)+\veps^2  c^{(1)}(\xi,t)+
\veps^4  c^{(2)}(\xi,t)+\cdots.}
\end{equation} 
{ 
Inserting the expansion~\ref{exprho} into Eq.~\eqref{newflux}, the following expansion for the flux hold true
\begin{equation}
    \veps J_\veps = J^{(0)} + \veps^2 J^{(1)}  + \veps^4 J^{(2)} + \cdots,
    \label{Jexpand}
\end{equation}
with 
\[
    J^{(k)} = -D\left(\pad{ c^{(k)}}{\xi} + c^{(k)}U'(\xi) \right), \quad k\geq 0.
\]
Using the expansion \ref{Jexpand}  in \ref{c_eps2}, the Eqs.~(\ref{eq_flux_eps}-\ref{eq_flux_eps_2}) become, to the various orders in $\veps$,
\begin{eqnarray}
    O(\varepsilon^{-1})  :&  \displaystyle \phantom{\quad\frac{\partial  c^{(k)}}{\partial t} + } \pad{J^{(0)}}{\xi} & = 0 ,
    \label{eq:eps0}\\
    O(\varepsilon^{k})  : &  \displaystyle \quad\frac{\partial  c^{(k)}}{\partial t} + \pad{J^{(k)}}{\xi}
     &  =  0, \quad k\ge 0.\label{eq:epsk}
\end{eqnarray}
where Eq.~\eqref{eq:eps0} states that  the lowest order flux $J^{(0)}$ is constant.\\
Condition $\lim_{x\to-\varepsilon}J(x,t) = 0$, which is a consequence of the singularity of the potential, gives:
\begin{equation}
\frac{\partial  c ^{(0)}}{\partial \xi}+ U'(\xi) c ^{(0)} = 0. 
\label{Jorder0} 
\end{equation} 
}
Therefore, it is possible to integrate the equation
\begin{equation}
\frac{1}{ c ^{(0)}}\frac{\partial  c^{(0)} }{\partial \xi} = - U'(\xi) 
\end{equation}
whose solution is
\begin{equation}
 c ^{(0)}(\xi,t) =  c ^{(0)}(L+1,t)\exp \left(-U(\xi) + U(L+1) \right) =  c ^{(0)}(L+1,t)\exp \left(-U(\xi)\right)
\end{equation}
where $U(L+1) = 0$ because the potential has effect only inside the bubble domain (that is $\Omega^\varepsilon_b$ for $V(x)$ and $[0,L+1]$ for $U(\xi)$), while it is constant outside, say $U(\xi)=0$ for $\xi \geq L+1$. 

The following assumption is made 
\[V_x(\varepsilon L) = 0,\] 
which is compatible {  with having} a smooth potential with compact support. Notice that, strictly speaking, LJ potential is not with compact support, however it is negligible small, together with its derivative, at distances which are {  sufficiently larger than its range $\varepsilon$.
}

Integrating \ref{pde1d} in the range of the potential, using the approximation $ c (x,t) =  c ^{(0)}(\xi,t), \, x \in \varepsilon[-1,L]$, the zero boundary condition \ref{bc1d} at $x=-\varepsilon $ and that $V_\varepsilon(\varepsilon L)=V'_\varepsilon(\varepsilon L)=0$, we obtain
\[
	\frac{d}{dt}\int_{- \varepsilon}^{\varepsilon L} c (x,t) \, dx +J_\varepsilon(\varepsilon L)  =  0,
\]
for which it follows
\begin{eqnarray}
	 \nonumber
 \frac{\partial c (\varepsilon L,t)}{\partial t} \; \varepsilon \int_{0}^{L+1}\exp\left(-U(\xi) \right) d \xi - D \frac{\partial  c (\varepsilon L,t)}{\partial x} &=&0
\end{eqnarray}
that represents a boundary condition at $x=\varepsilon L$.
 {Notice that, out of the range of the potential, the flux $J$ is given by just the diffusion term, thus Eq.~\eqref{eq:epsk} will not be used since for sufficiently small values of $\veps$ it is $c(x,t)\approx c^{(0)}(x,t)$, and Eq.~\eqref{eq:eps0} is used to provide an effective boundary condition for the diffusion equation.}

Summarizing, the \textit{multiscale model} can be stated as a limit model for $\varepsilon \to 0$ as follows:
\begin{eqnarray}\label{reduced1d}
\frac{\partial  c }{\partial t} &=& D\frac{\partial^2  c }{\partial x^2} \quad {\rm {in} }\, x \in [0,1]\\
\frac{\partial  c }{\partial x} &=& 0  \quad {\rm {at} }\, x = 1 \\ \label{BCt}
M\frac{\partial  c }{\partial t} &=& D\frac{\partial  c }{\partial x}  \quad {\rm {at} }\, x = 0
\end{eqnarray}
where
\begin{equation}
\label{expr_M}
M=\varepsilon\int_{0}^{L+1}\exp\left(-U(\xi)\right)d\xi.
\end{equation}

 {As $\veps\to 0$ one could let $L\to\infty$ in such a way that $L\veps\to 0$, so in the expression for $M$ the integral can be written from $0$ to $\infty$. 
In practice, $L$ is finite, and it is chosen $L=2$.  
}

\subsection{  Extension to more dimensions}
\begin{figure}[!ht]
	\centering
		\includegraphics[width=0.5\textwidth]{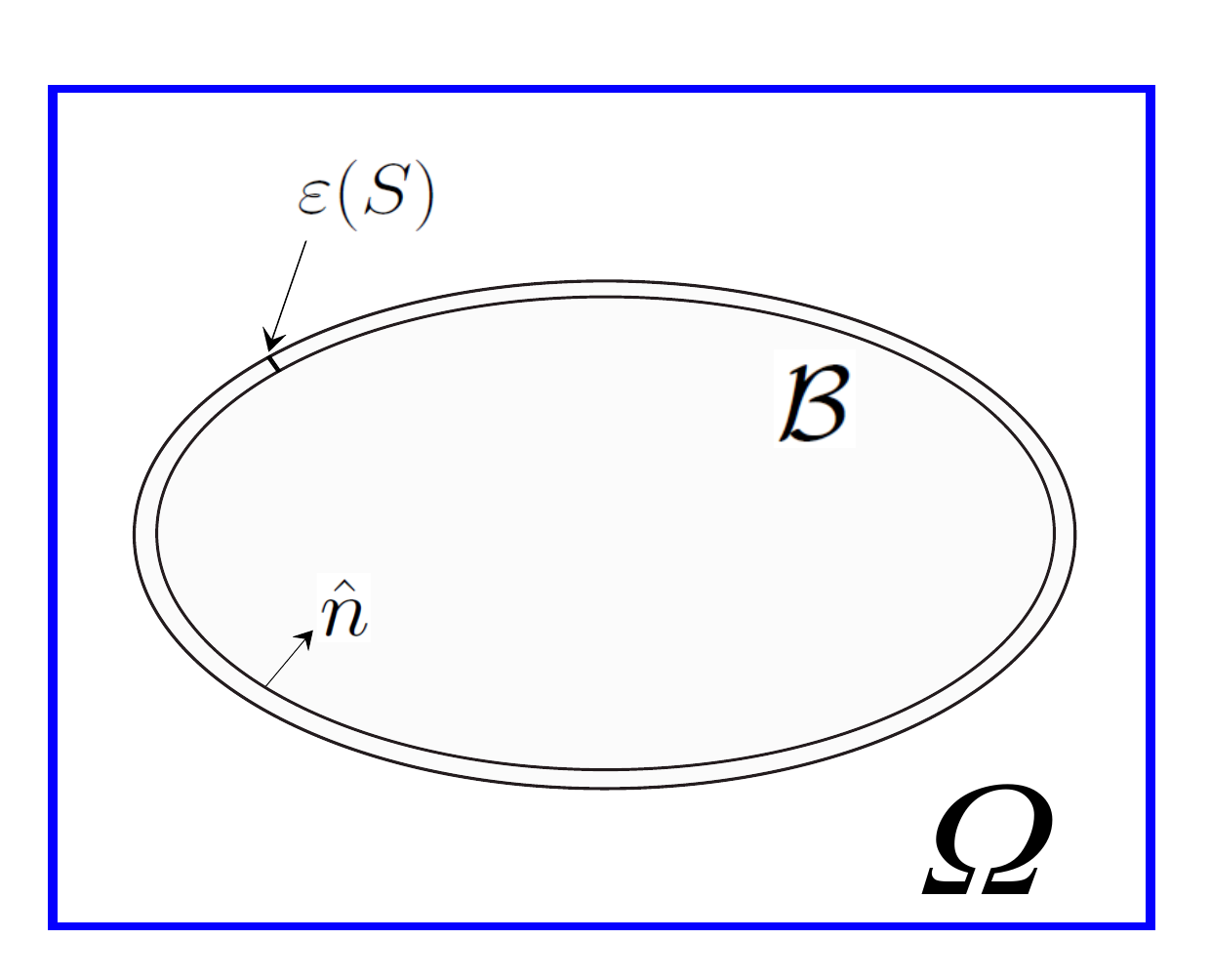}
	\caption{\textit{Scheme of the boundary condition in 2D where 
	$\hat{n}$ is the outgoing normal vector to the domain $\Omega$, { and $\varepsilon = \varepsilon(S)$ is the thickness parameter of the bubble surface $S$.  The internal line near the bubble surface  represents the boundary $\Gamma_B$.}}}
	\label{fig:bubble}
\end{figure}

{In this section the condition in Eq.~\eqref{BCt} is extended to 2D and 3D. The bubble $\mathcal{B}$ is defined by a level set function (see e.g. 
\cite{sussman1994level,sethian1999level,russo2000remark,osher2004level}) 
\begin{equation}
    \phi(\vec{x}) = r - \sqrt{x^2 + y^2},
\end{equation}
such that $\phi(\vec{x})< 0 $ inside the bubble and $\phi(\vec{x})> 0 $ outside. The unit normal $\hat{n}$ can be defined as 
\begin{equation}
    \hat{n}(\vec{x}) = \frac{\nabla \phi(\vec{x})}{|\nabla \phi(\vec{x})|} \qquad \forall \vec{x}\in\Gamma_B 
\end{equation}
where $\Gamma_B = \partial \mathcal{B}$.

Near the boundary, the flux can be decomposed into the normal and the tangential components, such that
\begin{equation}
    \vec{J} = J_n \hat{n} + J_\tau \hat{\tau}
\end{equation}
where $J_n $ and $J_\tau $ denote the components of $\vec{J}$ parallel and orthogonal to $\nabla V$, respectively.
The equation for the concentration $ c $ becomes:
\begin{equation}
	\frac{\partial  c }{\partial t} = -\left(\frac{\partial J_n}{\partial \hat{n}}+\frac{\partial J_\tau}{\partial \hat{\tau}}\right).
\end{equation}
Let $\vec{x}$ be a point in $\Gamma_B$ and, integrating along the normal direction, it becomes
\begin{eqnarray}
\label{eq_int_eps_x} \displaystyle 
	\frac{\partial}{\partial t}\int_{-\varepsilon}^{\varepsilon L} c (\vec{x}+r\hat{n})dr &= -J_n(\vec{x}-L\varepsilon\hat{n}) - \int_{-\varepsilon}^{\varepsilon L}\frac{\partial J_\tau}{\partial \tau}(\vec{x}+r\hat{n})dr
 \end{eqnarray}
since $J_n(\vec{x}+\varepsilon \hat{n}) = 0$. Now, we consider a change of variable from $\vec{x}$ to $\xi$ in the integration intervals and, making use of $\displaystyle J_\tau = -D\frac{\partial  c }{\partial \tau}$, the Eq.~\eqref{eq_int_eps_x} becomes
\begin{equation}
	\label{integrap}
	\varepsilon\frac{\partial}{\partial t}\int_{0}^{L+1} c (\vec{x}+\varepsilon \xi\hat{n})d\xi = -J_n(\vec{x}) + \varepsilon\int_{0}^{L+1}D\frac{\partial^2  c }{\partial \tau^2}(\vec{x}+\varepsilon\xi\hat{n})d\xi.
\end{equation}
Let us consider the term $J_n$:
\begin{equation}
	J_n = -D\left(\frac{\partial  c }{\partial \hat{n}}+\frac{1}{k_B T}\frac{\partial V}{\partial \hat{n}} c\right),
\end{equation}
since $\displaystyle V(x)/k_B T=U\left(\xi\right)$, it follows that 
\[\displaystyle \frac{\partial (V/k_B T)}{\partial \hat{n}} = \frac{1}{\varepsilon}U'(\xi)\] 
with $\xi = (r - R + \varepsilon)/\varepsilon$. For this reason, there is a dependence on $\xi $ of the solution inside the layer, $ c = c(\vec{x} + \varepsilon \xi \hat{n},t) = \tilde{c}(\xi,t)$, thus there is a common factor $1/\varepsilon$:
\begin{equation}
	J_n = -D\frac{1}{\varepsilon}\left(\frac{\partial  c }{\partial \xi}+ \frac{\partial U}{\partial \xi}c\right)
\end{equation}
where we omitted the $\tilde{\cdot}$. Following the same argument used in 1D to derive Eq.~\eqref{exprho}, the expansion in $\varepsilon$ of the solution, to the lowest order is
\begin{equation}
\label{bcbubble}
	\frac{\partial  c^{(0)}}{\partial \xi} + \frac{\partial U}{\partial \xi} c^{(0)} = 0.
\end{equation}
Integrating in $d\xi$, the solution of the Eq.\eqref{bcbubble} is
\begin{equation}
	c^{(0)}(\vec{x}+\varepsilon \xi\hat{n})= c^{(0)}(\vec{x})\exp\left(-U(\xi)\right).
\end{equation}
Using this expression we can compute the line density of entrapped ions:
\begin{equation}
	\mathcal{C}(\vec{x})\varepsilon\int_{0}^{L+1} c (\vec{x}+\varepsilon \xi\hat{n})d\xi \simeq M c (\vec{x})
\end{equation}
with $\displaystyle M= \varepsilon \int_{0}^{1}\exp\left(-U(\xi)\right)d\xi$ and where we assumed $U(L+1)=U'(L+1)=0$.

At the end, in 2D (and also in 3D), the system and the new boundary condition for the concentration becomes:
\begin{eqnarray}	\label{system_multiscale}
	\displaystyle \frac{\partial c}{\partial t} &=& D \Delta c \quad \rm{in} \> \Omega \\ \nonumber
	\displaystyle \frac{\partial c}{\partial n} &= &0 \quad \rm{ on } \, \partial \Omega\backslash \Sigma_B \\ \label{bc_multiscale}
	\displaystyle M\frac{\partial  c }{\partial t} &=& MD\Delta_\perp  c +D\frac{\partial  c }{\partial n} \quad \rm{ on }\, \Sigma_B
\end{eqnarray}
because $J_n = -D\partial  c /\partial n$ just out of the range of the potential.

Notice that in 2D the Laplace-Beltrami operator reduces to the second derivative with respect to the arclength of the boundary:\footnote{Here by $\partial/\partial \tau$ 
we denote the derivative on $\Gamma_B$, i.e.\ the derivative along the arclength that parametrizes the curve, likewise
$\partial^2/\partial \tau^2$ denotes the second derivative along $\Gamma_B$, not the second derivative along the tangent direction. See Fig.\ \ref{fig:bubble} and Appendix~\ref{A:LB}. }
\begin{equation}
    \Delta_\perp c = \frac{\partial^2 c}{\partial \tau^2}.
    \label{eq:BL2D}
\end{equation}

}

\subsection{Validation of the model and results}
\label{section_validation_multiscale}
In this section we validate the model in 1D and 2D.

First, a 1D domain is considered as intermediate step. Then, a more realistic 2D domain is considered, where a circular hole, at the center of a squared domain, describes the attractive bubble surface.

\subsubsection{Validation 1D}
The initial domain in 1D is $\Omega^\varepsilon= [- \varepsilon,1] =\Omega^\varepsilon_b \cup \Omega^\varepsilon_f ,\,\Omega^\varepsilon_b = [- \varepsilon,2\varepsilon],\, \Omega^\varepsilon_f = [2\varepsilon,1]$, 
and the space derivatives are discretized with Finite-Difference (FD) schemes, see Chapter ?? for more details. In this case, we refer to the system (\ref{pde1d}-\ref{expr_V_LJ}). 

For the multiscale model, the domain is $\Omega^0 = [0,1]$, and the equations are in (\ref{reduced1d}-\ref{expr_M}). 
Internal equations are discretized by finite differences by the classical three points approximation of the second derivative, while the Eq.~\eqref{BCt} is discretized to second order as 
\begin{equation} 
\label{eq_BCt_discr}
    M\frac{\partial }{\partial t} \frac{c_0 + c_1}{2} = D \frac{c_1 - c_0}{h}. 
\end{equation}

\begin{figure}[tb]
	\begin{minipage}
		{.44\textwidth}
		\centering
		\includegraphics[width=\textwidth]{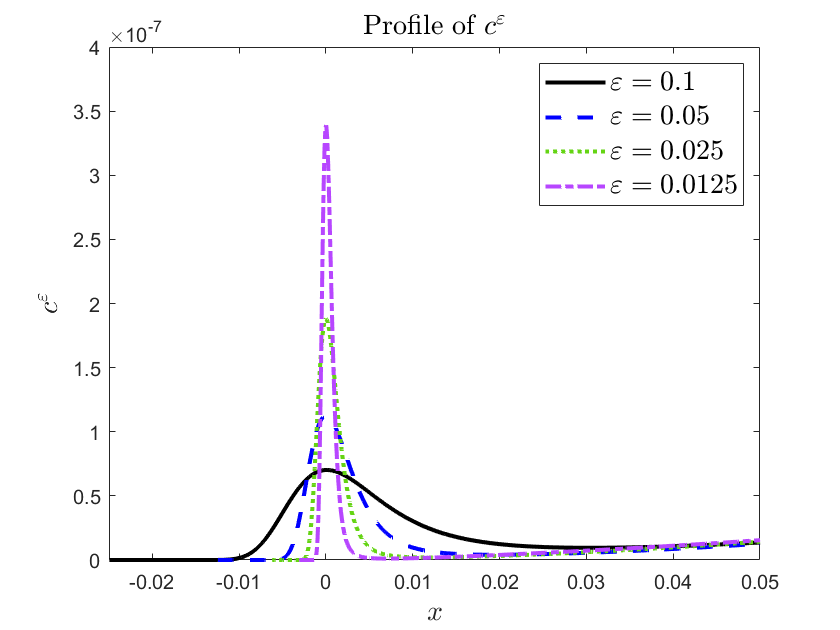}
	\end{minipage}
	\begin{minipage}
		{.44\textwidth}
		\includegraphics[width=\textwidth]{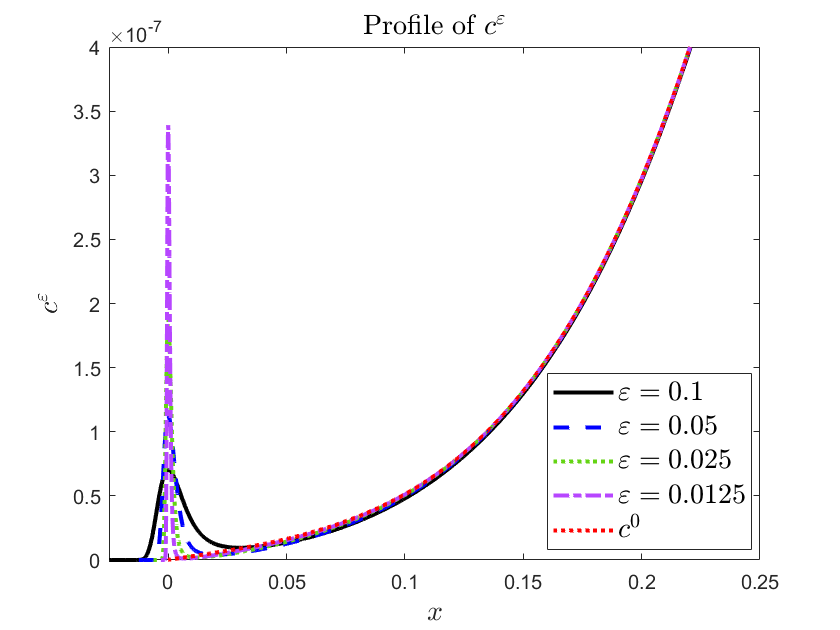}
	\end{minipage}
	\caption{\textit{Here there are two aspects of the same plot. On the left it is shown the effect of the attractive potential that becomes higher decreasing the value of $\varepsilon$. On the right there is the comparison between $c^\varepsilon$, and the dashed red line that represents the solution of the multiscale model, $c^0$.}} 
\label{1D_full_complete}
\end{figure}

Let $c^\varepsilon$ and $c^0$ be the solutions of the full and multiscale models, respectively, and their comparison is showed in Fig.~\ref{1D_full_complete}. 

\subsubsection{Validation 2D in presence of a circular bubble}
\label{sec:validation2D}
The arbitrary boundary of the bubble is described by a level set function. The set of equations is then solved by a finite difference scheme on a regular Cartesian mesh, which makes use of regular grid points (inside the domain $\Omega^0$) and ghost points (neighbors of regular grid points which are out of $\Omega^0$). The equations on the regular points are obtained by finite difference discretization of the Laplacian operator, while the equations on each ghost point are given by interpolating the numerical solution near the ghost point and imposing the boundary condition on the interpolant at the orthogonal projection of the point on $\partial \Omega^0$ (see \cite{COCO2013464,COCO2018299} for more details).

\begin{figure}[htb]
	\centering	
	\begin{minipage}
		{.4\textwidth}
		\centering
		\includegraphics[width=\textwidth]{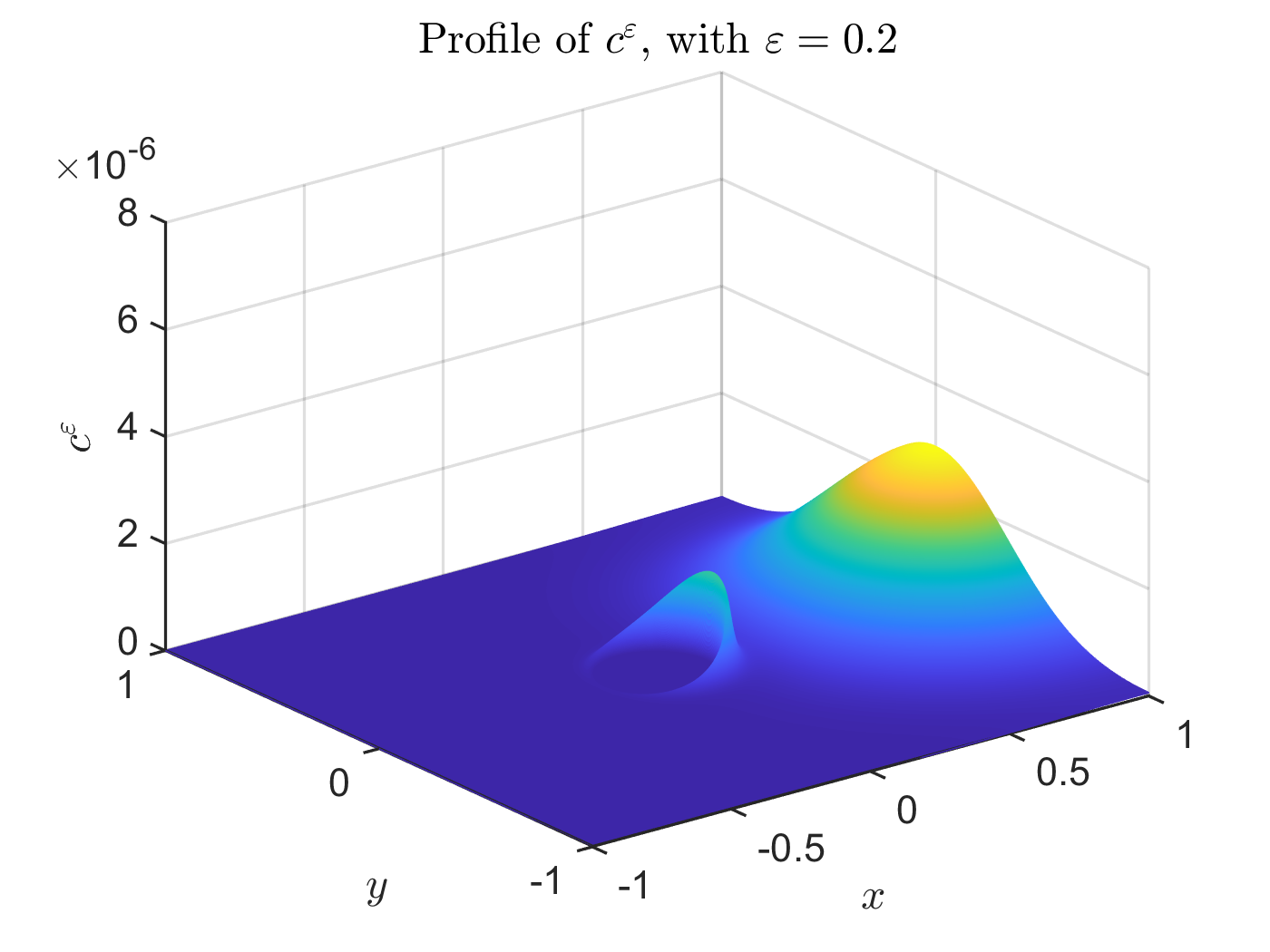}
	\end{minipage}
	\begin{minipage}
		{.4\textwidth}
		\centering
		\includegraphics[width=\textwidth]{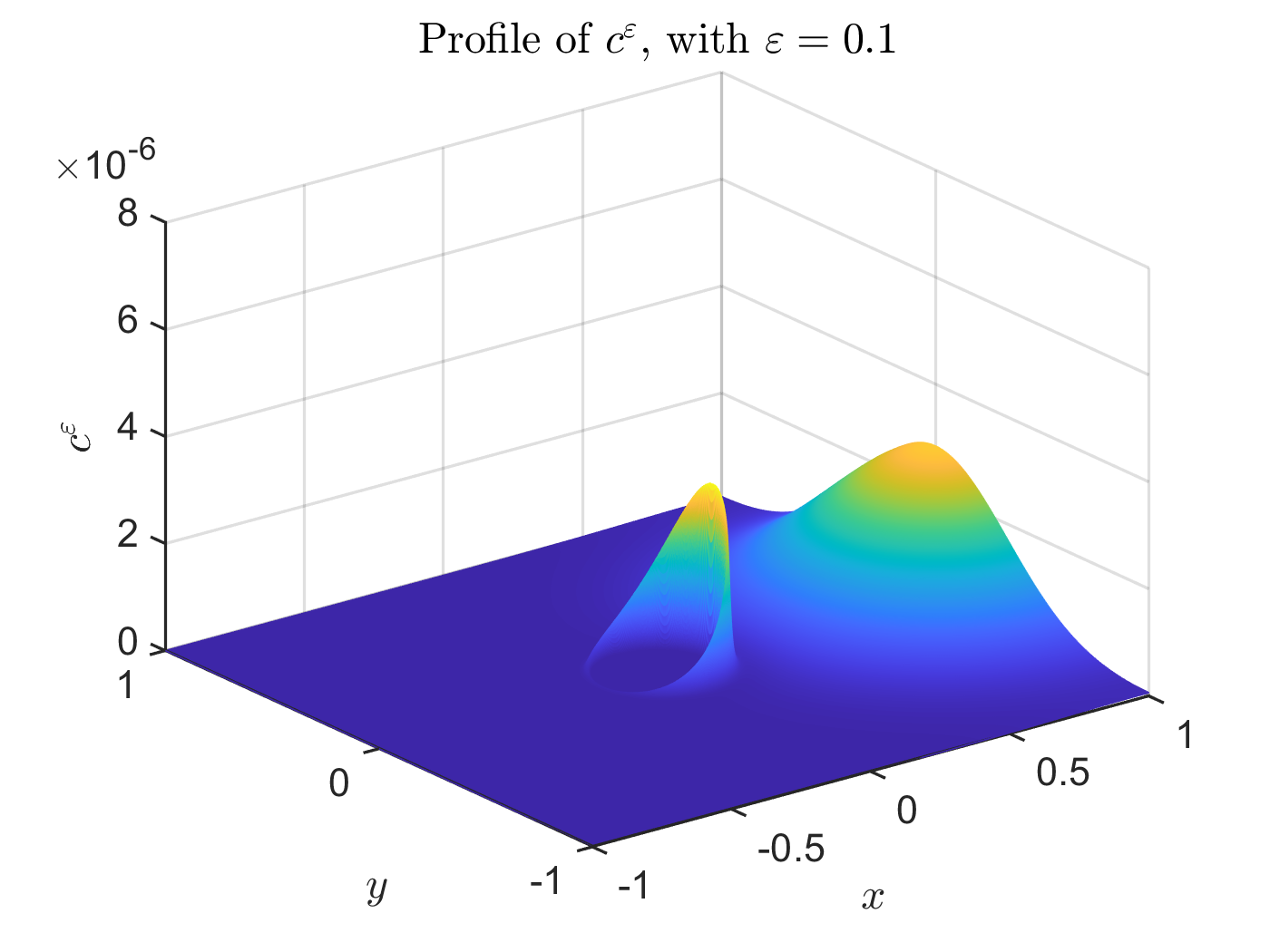}
	\end{minipage}
	\begin{minipage}
		{.4\textwidth}
		\centering
		\includegraphics[width=\textwidth]{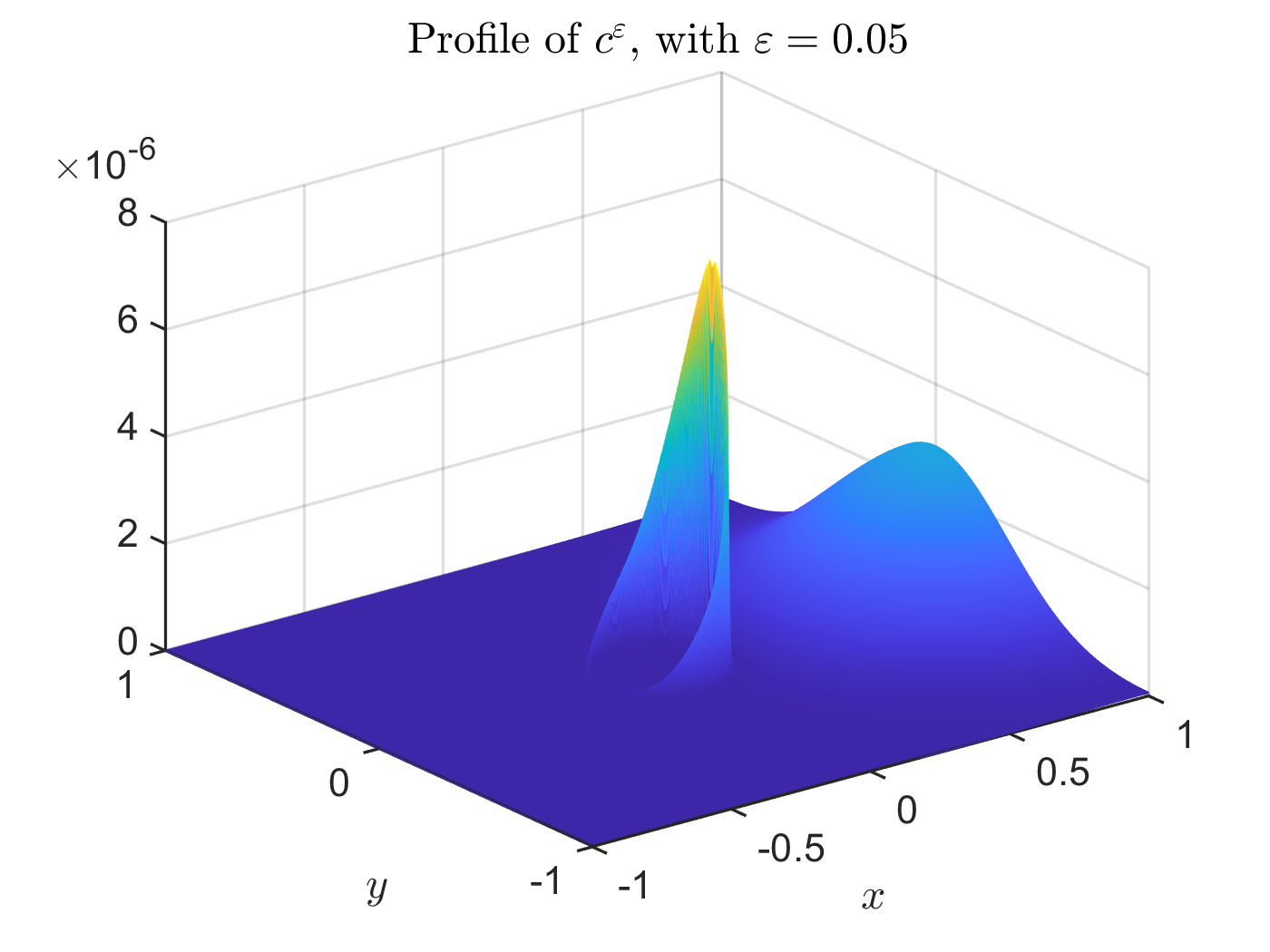}
	\end{minipage}
	\begin{minipage}
		{.4\textwidth}
		\centering
		\includegraphics[width=\textwidth]{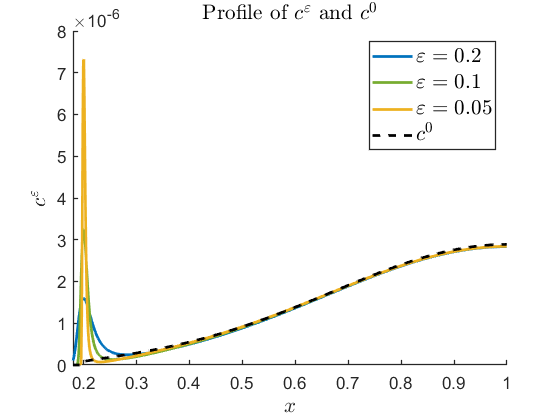}
	\end{minipage}
	\caption{\textit{The effect of the attractive potential in a 2D domain. Here the bubble has a reasonable circular shape. As before, in the first three panels (a,b,c) there are the solution $c^\varepsilon$ for three different values of $\varepsilon$. In panel (d) there is a 1D plot, that is a section of the previous ones for $y=0$ with an additional line. The black dashed line represents the solution of the multiscale model, $c^0$.
 }}
\end{figure}

\section{Coupling conditions for cardiovascular system}
\label{sec:cardiovascular}
Talking about the circulatory system, it is very important to introduce the concept of geometric multiscale. 
In this case, it is related to the coupling of different dimensions and scales, in the same system of equations. Considering the relationships between different scales and dimensions, it is possible to gain a more comprehensive understanding of the system.

\begin{figure}
\centering
 \begin{overpic}[abs,width=0.6\textwidth,unit=1mm,scale=.25]{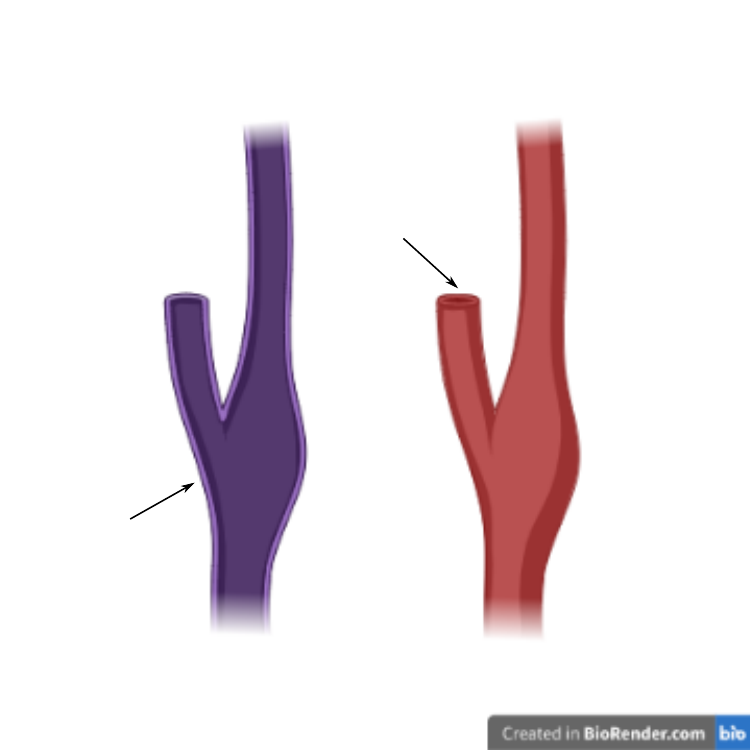}
\put(10,22){{$\Sigma^t$}}
\put(18,65){{$\Omega^t_f$}}
\put(48,65){{$\Omega_s$}}
\put(40,55){{$\Sigma$}}
\end{overpic}
	\caption{\textit{An example of a chopped domain. On the left, the purple object represents the time dependent fluid domain, while the red one represents the structure domain. }}
	\label{fig_chopped_domain}
\end{figure}

The cardiovascular system is very complex for different reasons. Apart from the geometric complexity, it is very important to mention the difficulty in obtaining data that are used to deduce meaningful initial and boundary conditions, when an artificial 'chopped' domain is considered for computational purposes (see, e.g. Fig.~\ref{fig_chopped_domain}). Local artificial domains may have global consequences, and to take into consideration these properties, suitable boundary conditions need to be applied between the free surfaces of the artificial domain.  

\subsection{The 3D model}\index{The 3D model}
\begin{figure}
\centerline{\includegraphics[width=0.5\textwidth]{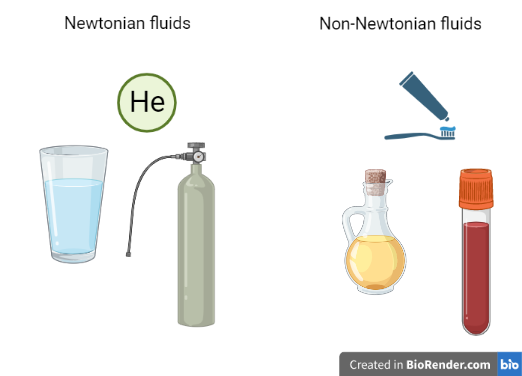}}
	\caption{\textit{Classification of Newtonian and Non-Newtonian fluids}}
	\label{fig_new_nonnew}
\end{figure}
In order to obtain a realistic description of a vessel flow, it is convenient to refer to a 3D Navier-Stokes formulation for incompressible homogeneous Newtonian fluids. Since the diameter of blood cells is about $10^{-5}m$, while the smallest diameter of arteries and veins is of the order of $10^{-3}m$, linear equations govern the relationships between internal forces and velocity gradient, {and blood can be treated as a  called {\it Newtonian fluid\/} } (see, e.g. \cite{perktold1986numerical,formaggia2009multiscale}). In Fig.~\ref{fig_new_nonnew} we see few examples of Newtonian (water, Helium and Oxygen) and Non-Newtonian fluids (tooth paste, oil, blood). As it is shown in that figure, the blood belongs to the Non-Newtonian fluids, but with the approximations mentioned above, it is possible to consider it as a Newtonian fluid.

An \textbf{Eulerian} formulation is more convenient for the fluid equations and the frame of reference is denoted by $\Omega_{\rm f} \subset \mathbb{R}^3$  (see Fig.~\ref{fig_chopped_domain} (left panel)). While, for deformable solid structure (the arteries) a \textbf{Lagrangian} framework is more suitable  (see Fig.~\ref{fig_chopped_domain} (right panel)) and the reference domain is $\tilde{\Omega}_{\rm s} \subset \mathbb{R}^3$. Therefore, there exists a map $\mathcal{L}:\tilde{\Omega}_{\rm s} \to {\Omega}_{\rm s}$, where $\Omega_{\rm s}$ is the spatial domain.
$\Sigma$ will be the interface between the fluid and the structure domains.

Now it is possible to define the system for the velocity field $\vec{u}$, the fluid pressure $p$ and the density $\rho$, both for the fluid and the structure parts, such that, 
\begin{equation}
\label{sys_3D_mass_momentum}
\partial_t {\bf V}_{\rm l} + M_{\rm l}({\bf V}_{\rm l}) \, \nabla \cdot {\bf V}_{\rm l} = 0, \quad M_{\rm l}({\bf V}_{\rm l}) \in \mathbb R^{3,3}, \quad {\rm l = f,s}
\end{equation}
where $l = f$ stands for the fluid system and variables and $l = s$ for the structure ones. The vectors of unknowns are
\begin{equation*}
{\bf V}_{\rm l}^{\text{T}} = (\rho_{\rm l}, u_{\rm l}, p_{\rm l}), \quad {\rm l = f,s}.
\end{equation*}
If the matrix $M_{\rm l}, \, {\rm l = f,s}$ admits real eigenvalues, and it is diagonalizable (i.e., the eigenvectors form a bases) the system is called hyperbolic, and it is well-posed.

\begin{example}
Let us consider Euler's system, in the case of ideal gases. The equation for the conservation of mass is the following
\begin{equation*}
\frac{\partial \rho}{\partial t} + u_i \partial_i \rho + \rho \partial_i u_i = 0 \quad\Rightarrow\quad \frac{D\rho}{Dt} + \rho \nabla \cdot \vec u = 0,
\end{equation*}
the one for the conservation of the momentum reads
\begin{equation*}
\rho \frac{D u}{Dt} + \nabla p = 0, 
\end{equation*}
while the relation for the energy can be reduced as
\begin{equation}
\label{eq:sound speed}
\frac{D p}{Dt} = \biggl(\frac{\partial p}{\partial \rho}\biggr)_S \frac{D \rho}{Dt} = a^2 (-\rho \nabla \cdot \vec u), \quad a^2 :=\biggl(\frac{\partial p}{\partial \rho}\biggr)_S.
\end{equation}
In this way we obtain a system for the variables $\rho$, $p$ and $\vec u$:
\begin{equation}
\label{sist:fondamentale-2}
\left \lbrace
\begin{aligned}
& \frac{D\rho}{Dt} + \rho \nabla \vec u = 0  \\
& \frac{D u}{Dt} + \frac{1}{\rho} \nabla p = 0  \\
& \frac{D p}{Dt} + a^2 \rho \nabla \cdot \vec u = 0\\
\end{aligned}
\right.
\end{equation}

For simplicity we consider the system \eqref{sist:fondamentale-2}, in 1D, that becomes
\begin{equation}
\label{sist:start}
\left \lbrace
\begin{aligned}
\partial_t \rho &+ u \,\partial_x \rho + \rho\, \partial_x u &= 0 \\
\partial_t u &+ u \,\partial_x u + \frac{1}{\rho}\, \partial_x p &= 0 \\
\partial_t p &+ u\, \partial_x p + a^2 \rho  \, \partial_x u &= 0 \\
\end{aligned}
\right.
\end{equation}
Now we define the following vector
\begin{equation*}
V^{\text{T}} = (\rho, u, p)
\end{equation*}
and the system \eqref{sist:start} can be rewritten as
\begin{equation}
\label{def:matrice}
\partial_t V + M(V) \, \partial_x V = 0, \quad M(V) = \begin{pmatrix} u & \rho & 0 \\ 0 &u&1/\rho \\ 0 & a^2 \rho & u \end{pmatrix}
\end{equation}

Now we calculate the three eigenvalues of the matrix $M$:
\begin{equation}
\label{eq:autovalori dx}
\text{det}(\lambda \mathbb I - M) = (\lambda - u) \bigl [ (\lambda - u)^2 - a^2 \bigr ] = 0 
\end{equation}
that are $\lambda_0 = u$ e $\lambda_{\pm} = u \pm a$, and they belong to $\mathbb R$ if $a \in \mathbb R$. 

If $a \in \mathbb C$ the system is not hyperbolic!
\end{example}

\subsection{The 1D and 0D model}\index{The 1D and 0D model}

Modern supercomputers struggle with the complexity of creating a complete 3D model of the cardiovascular system. The computational process for such a model takes hours or even days, whereas 1D or 0D models achieve good accuracy in a much shorter time frame. Interestingly, in certain situations, the accuracy provided by the 3D models exceeds the requirements of the problem at hand.

The focus should be on obtaining a detailed depiction of blood flow in specific critical areas where a pathology is identified or suspected. For this purpose, a highly efficient approach involves utilizing a comprehensive 3D model solely in the designated region, while employing simplified 1D or 0D models elsewhere \cite{quarteroni2016geometric,quarteroni2017cardiovascular}.
The 1D-branches are then connected to the rest of the circulatory system by 0D equations, i.e. by a system governed by ODE's.

\begin{figure}
\centerline{\includegraphics[width=0.5\textwidth]{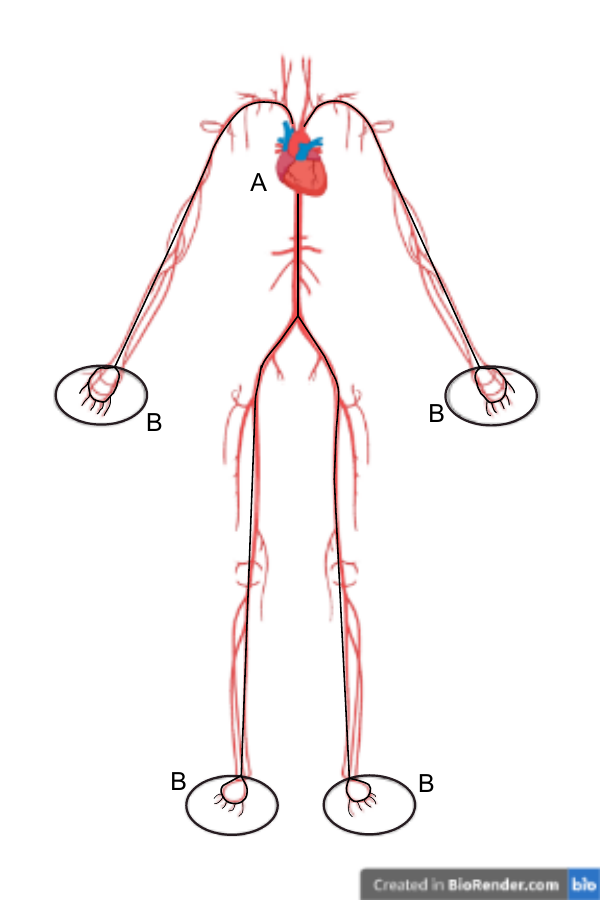}}
	\caption{\textit{Systemic arterial tree consisting of larger arteries originates at the heart (A) and
terminates at the points marked with (B) and structured trees representing
smaller arteries originate at these terminals and provide main
tree with outflow boundary conditions.}}
	\label{fig_1D_trees}
\end{figure}

At the end, the entire domain will be a combination of several domains of 3D, 1D and 0D models. What is important is the continuity relations at the junctions between them. 
The remaining domain is represented by attaching a structured tree, to simulate the spreading of certain tissues, see Fig.~\ref{fig_1D_trees}.

To deduce a 1D model, some assumptions are needed:
\begin{enumerate}
    \item the axis $z$ of the cylinder is assumed rectilinear and the component $u_z$ of the velocity is the dominant one
    \item all the cross sections of the cylinders are circles $\mathcal{S}(t,z)$ with radius $R(t,z)$, and the pressure is considered constant over each section $\mathcal{S}(t,z)$
    \item the membranes have constant thickness
\end{enumerate}




To derive the 1D model, it is enough to integrate the Eqs.~(\ref{sys_3D_mass_momentum}) over the section $\mathcal{S}$ in the $z$-direction. 
An additional order reduction is obtained integrating the system along the axis direction, obtaining the 0D model. In this kind of models, only the time dependence survives, and they are common to describe phenomena in capillaries. 


\subsection{Coupling of 3D-1D, 3D-0D and 1D-0D models}\index{Coupling of 3D-1D, 3D-0D and 1D-0D models}
To accurately describe the entire cardiovascular system, it is crucial to select the appropriate spatial dimensions that adequately represent each component, depending on the desired level of accuracy. In certain regions, where precision is of utmost importance, a 3D model is employed. However, in other areas, computational complexity can be reduced by utilizing simplified models such as 1D or 0D.

\begin{figure}
\centering
\begin{overpic}[abs,width=0.35\textwidth,unit=1mm,scale=.25]{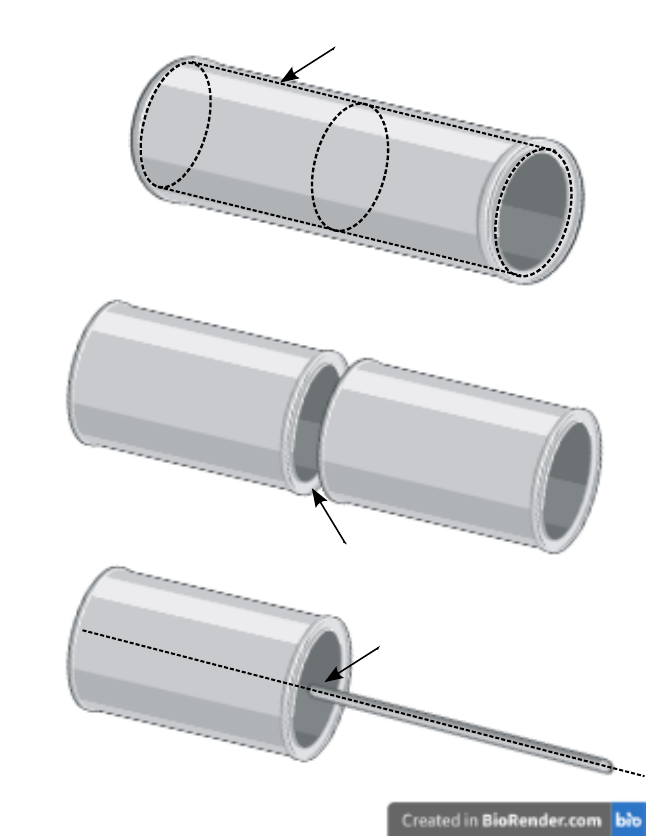}
\put(-7,56.5){(1)}
\put(22,51.5){$\Sigma$}
\put(17,43){$\Omega$}
\put(-7,38){(2)}
\put(12,27){$\Omega_{\rm f}$}
\put(28,23){$\Omega_{\rm s}$}
\put(-7,18){(3)}
\put(23,16){$\Gamma = \{\Gamma_{\rm f},\Gamma_{\rm s}\}$}
\put(25,11){{$z=0$}}
\put(42,6.5){{$z=L$}}
\end{overpic}
\caption{\textit{Scheme of domain splitting (the subscript ${\rm f}$ stands for fluid and ${\rm s}$ for solid): (1) $\Omega$ is the initial considered domain and $\Sigma$ the related boundary. (2) $\Omega_{\rm f}\cup \Omega_{\rm s}$ is a partition of $\Omega$ and $\Gamma$ is the common interface. (3) The 3D domain $\Omega_{\rm s}$ is replaced with a 1D domain $[0,L]$ in $z$-direction. }}
	\label{fig_domain_splitting}
\end{figure}

When combining both 3D and 1D (or 0D) models, a combination of these models is utilized in specific regions, requiring an approximation to ensure continuity of the quantities involved. To facilitate this process, the first step involves a technique called \textit{domain splitting}. It implies the division of a cylindrical vessel into two subdomains, as depicted in the accompanying Fig.~\ref{fig_domain_splitting}. Subsequently, the 3D-1D (or 3D-0D) model can be formulated by incorporating appropriate interface conditions.


Regarding the 3D-0D coupling, there are three properties that mathematical modeling would face: (i) the non linear terms are neglected; (ii) the boundary conditions will be integrated in the equations; (iii) domain splitting (see Fig.~\ref{fig_domain_splitting}) is required, and some continuity property can be lost (the subscript ${\rm f}$ stands for fluid and ${\rm s}$ for solid). 

At the end, we show the strategy of the 1D-0D coupling. Suppose there is a 1D model 
\begin{equation}
\label{eq_1D_0D_zdim}
    \mathcal{M}_{1D}({\bf V}_{1D},\Omega) = 0
\end{equation}
where ${\bf V}_{1D}$ is the vector variables and $\Omega = [-L,L]$ is the domain. The problem  can be split in two sub-problems $\mathcal{M}^1_{1D}$ and $\mathcal{M}^2_{1D}$ and the domain in two subdomains $\Omega_1 = [-L,0]$ and $\Omega_2 = [0,L]$: 
\begin{eqnarray}
\label{eq_1D_1D_1}
    \mathcal{M}^1_{1D}({\bf V}^1_{1D},[-L,0]) = & 0, \quad {\rm in } \, \Omega_1 \\ \label{eq_1D_1D_2}
    \mathcal{M}^2_{1D}({\bf V}^2_{1D},[0,L]) = & 0, \quad {\rm in } \, \Omega_2 \\
    \label{eq_1D_1_con}
    {\bf V}^1_{1D}(z=0) = & {\bf V}^2_{1D}(z=0)
\end{eqnarray}
where Eqs.~\eqref{eq_1D_1D_1},\eqref{eq_1D_1D_2} are the equations for the domains $\Omega_1$ and $\Omega_2$, respectively, while Eqs.~\eqref{eq_1D_1_con} is the continuity conditions at the intersection between $\Omega_1$ and $\Omega_2$, $z=0$. 
The next step is to replace the 1D problem in $\Omega_2$, with a suitable 0D model, with appropriate boundary conditions. To do this, average quantities are needed in 0D, leading to the following problem
\begin{eqnarray}
      \mathcal{M}^1_{1D}({\bf V}^1_{1D},[-L,0]) = & 0, \quad {\rm in } \, \Omega_1 \\
 {\bf V}^1_{1D}(z=0) = & {\bf V}_{0D} \\
    \mathcal{M}_{0D}({\bf V}_{0D}) = & 0.
\end{eqnarray}

\section{Coupling conditions for shallow water equations at channel junctions}
\label{sec:SW_channels}
The shallow water equations is a hyperbolic system that describes surface water flows, such as in artificial channels, rivers, lakes and oceans. The system consists of a mass and a momentum balance equation, and, in the most common scenario, it is considered in  a two dimensional domain. 

This section is devoted to the description of an efficient way to consider one dimensional problems. In Fig.~\ref{fig:junction} we see the schematic representation of the intersection of three channels. As a reference example, we consider the first channel the one that is parallel to the $x$-axis, and the other two create two angles with the first one. The figure shows how to represent with a 1D domain a 2D problem.

As we said in the Introduction, the tricky part is to define the boundary condition of each channel at the junction, and to couple them with the SW equations of the three channels. The main idea is to see the junctions as in the Fig.~\ref{fig:junction}, focusing at the intersections with a more detailed 2D domain.

Choosing the coupling conditions it can be seen as a Riemann problem defined by imposing conservation of mass and continuity of the water height through junctions. To explain better the concept, we can say that the following problem is solved at the junctions:
\begin{equation}
\label{eq_RS_junctions}
\left\{ \begin{array}{c} \displaystyle \partial_t u_l + \partial_x f(u_l) = 0, \quad t \in \mathbb{R}^+, \quad l =1,2,3  \\ 
    u_l(0,x) = \bar{u}_l, \quad x\in \mathbb{R}^+, \, u_l \in \Omega, \quad l =1,2,3.
    \end{array}  \right.
\end{equation}
where $\bar{u}_l, l =1,2,3$ are the constant states in $\Omega$. $u_l$ has two components, i.e., $u_l = (u_{l,L},u_{l,R}), l = 1,2,3,$ that denote the densities of the conserved quantities in the $l-th$ tube. 
\subsection{Mathematical model}
As we said before, the main strategy is to consider 1D models along the straight channels, and 2D models at the junctions, in order to obtain detailed boundary conditions. 

The shallow water equations in 2D read:
\begin{equation}
\label{eq:junction_2D}
    \partial_t {\bf Q}_{2D} + \partial_x {\bf F}_{2D}({\bf Q}_{2D}) + \partial_y {\bf G}_{2D} ({\bf Q}_{2D}) = {\bf S}({\bf Q}_{2D})
\end{equation}
where ${\bf Q}_{2D} = [h,hu,hv]^T$ is the vector of conserved variables, ${\bf F}_{2D}({\bf Q}_{2D}) = [hu, hu^2 + \frac{1}{2} gh^2, huv]^T$ is the flux in $x$-direction and ${\bf G}_{2D} ({\bf Q}_{2D}) = [hv,hvu,hv^2 + \frac{1}{2}g h^2]^T$ is the flux in $y$-direction, $g$ is the gravity and  ${\bf S}_{2D}({\bf Q}_{2D}) = [0, -ghS_x,-ghS_y]^T$ is the variation of the bottom topography, with $S_x, S_y$ known functions.

The same model, in 1D, reads
\begin{equation}
\label{eq:junction_1D}
    \partial_t {\bf Q}_{1D} + \partial_s {\bf F}_{1D}({\bf Q}_{1D}) = {\bf S}_{1D}({\bf Q}_{1D})
\end{equation}
where the second component of the velocity $v$ is dropped, the $y$-direction the flux {\bf G} disappears and $u_{1D}$ is the velocity along the $s$-direction.

It is possible to combine the 1D system, in Eq.~\eqref{eq:junction_1D}, defined in a generic direction $s$, with the 2D system defined in Eq.~\eqref{eq:junction_2D}. We start integrating Eq.~\eqref{eq:junction_2D}, that becomes 
\begin{eqnarray}
   \frac{\partial }{\partial t}\int_V {\bf Q}_{2D} dV + \int_\Omega \left(\cos (\theta) {\bf F}_{2D} ({\bf Q}_{2D})  + \sin(\theta) {\bf G}_{2D}({\bf Q}_{2D})  \right) d\Omega = 
   {\bf 0}.
\end{eqnarray}
where $[\cos \theta, \sin \theta]$ are the two components of the outward unit normal vector, such that, ${\bf n} = [n_1,n_2] = [\cos \theta, \sin \theta]$, and $V$ is the control volume in the Cartesian plane. 

If ${\bf T} = {\bf T}(\theta)$ is a rotation matrix, and ${\bf T}^{-1}(\theta)$ its inverse, such that
\begin{equation}
\label{eq:T_theta}
T=\left[\begin{array}{ccc} 
1 & 0 & 0\\
0 & \cos(\theta) & \sin(\theta)\\ 
0 & -\sin(\theta) & \cos(\theta)\\
\end{array}\right], \qquad T^{-1}=\left[\begin{array}{ccc} 
1 & 0 & 0\\
0 & \cos(\theta) & -\sin(\theta)\\ 
0 & \sin(\theta) & \cos(\theta)\\
\end{array}\right] 
\end{equation}
it is easy to prove that 
\begin{equation}
\label{eq_TFT}
    {\bf H} = {\bf n}\cdot [{\bf F}({\bf Q}),{\bf G}({\bf Q})] = \cos (\theta) {\bf F}({\bf Q})  + \sin(\theta) {\bf G}({\bf Q}) = {\bf T}^{-1}{\bf F} ({\bf T}({\bf Q}))
\end{equation}
for all vectors ${\bf Q}$ and all normal directions of $\Omega$

The final steps to describe are the following:
\begin{itemize}
    \item 1D system (Eq.~\eqref{eq:junction_1D}) where the direction $s$ is the normal direction, with rotated Riemann initial condition
    \item calculate the angle $\theta_e$ between the outward unit normal to the edge $e$ and the reference $x$ direction, being positive in the counter-clock direction
    \item calculate the corresponding rotation matrix ${\bf T_e}$ and rotate left and right initial data
    \item solve the 1D Riemann problem and compute an approximation of the flux in the edge $e$
    \item rotate back the flux (see Eq.~\eqref{eq_TFT})
\end{itemize}

\begin{figure}[tb]		
\centering
\begin{overpic}[abs,width=0.55\textwidth,unit=1mm,scale=.25]{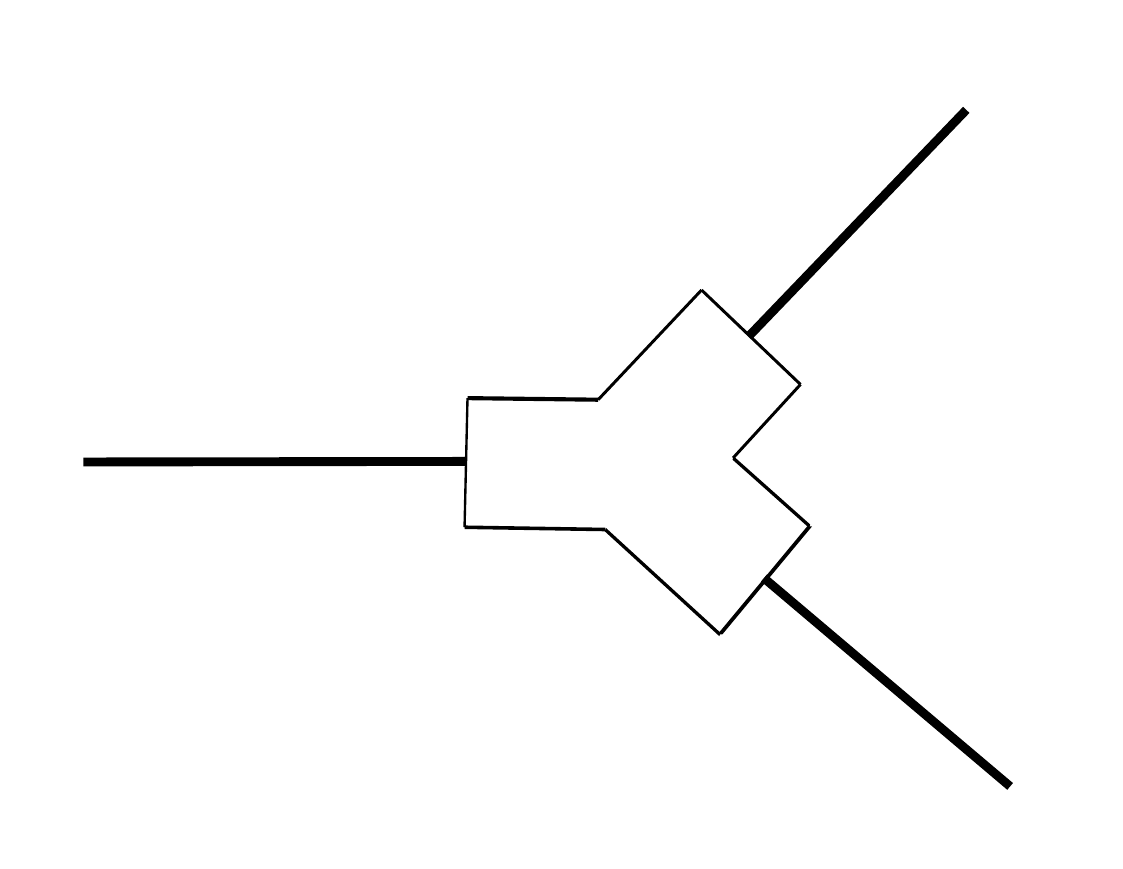}
\put(18.5,30){1D}
\put(52,45){1D}
\put(57,16){1D}
\put(38,26){2D}
\end{overpic}
\caption{\textit{Coupling 2D and 1D domain.}}	
	\label{fig:junction_2}
\end{figure}


\subsection{Methodology}
The basic idea of the first method is to consider a 1D formulation for the straight channels and a single 2D element at the junctions, as shown in Fig.~\ref{fig:junction_2}. It is important to explain how to switch from 1D to 2D model and vice-versa. 

The 2D flux is calculated by solving the rotate 1D Riemann problem \eqref{eq_RP_1D} in local coordinates, that reads as follows:
\begin{equation}
\label{eq_RP_1D}
\left\{ \begin{array}{c}  \displaystyle \partial_t {\bf Q}_{1D} + \partial_s {\bf F}_{1D}({\bf Q}_{1D}) = {\bf 0}. \\ 
      \displaystyle {\bf Q}_{1D}(s,0) = \left\{ \begin{array}{c}  \displaystyle  {\bf Q}_{1D,L} = {\bf T}_e({\bf Q}_{2D,L}) \quad {\rm if } \, s < 0 \\ 
    \displaystyle  {\bf Q}_{1D,R} = {\bf T}_e({\bf Q}_{2D,R}) \quad {\rm if } \, s > 0
\end{array}  \right.
    \end{array}  \right.
\end{equation}
where, the first line is the PDEs in normal direction, while the other condition are the rotated initial conditions.

Near to the junctions, a 2D system of equations is required. Referring to Fig.~\ref{fig:junction}(a) the result is a 1D system, and 
it is possible to obtain a non-zero transversal components for the 1D element, as follows
\begin{equation}
   \displaystyle  u_{1D} = {\rm sign}\left( u_{1D} \right)\sqrt{ u_{2D}^2 + v_{2D}^2 }.
\end{equation}

\section{Appendix}
\subsection{Shallow water equations}
\label{section_sw}
Shallow water theory is an approximation where the depth of the water is assumed to be sufficiently small compared to some significant length, like the radius of curvature of the water surface. This theory applies to many occurrences in nature, such as flows in open channels, like roll waves, and flood waves in rivers. It is commonly used in hydraulics by engineers dealing with these applications. 
{It can be shown that under such conditions, the pressure in the fluid region is essentially the hydrostatic pressure. 

The formal way to derive SW equations, starting from the exact hydrodynamical system, is obtained as the lowest order approximation in the asymptotic expansion procedure. The goal is to construct and justify the asymptotic expansion of the solution in powers of the characteristic quantity defined as the ratio between the typical length in $y$-direction and the one in $x$- (and $z$-) direction. In \cite{stoker1992water} Section 2.4, there is an explicit derivation of the $0^{th}$ and $1^{st}$ orders of the asymptotic expansion, and, after some integration, one is able to re-obtain the hydrostatic pressure relation that is the analogue of the assumption made for the SW equations.  

Here we provide a simple derivation of the shallow water equation under the assumption of irrotational flow and {\it hydrostatic approximation}. We restrict to consider a two dimensional problem in the $x$-$y$ vertical plane.} 

The $x$-axis is represented by the undisturbed free surface of the water, and the $y$-axis is taken vertically upwards. Thus, the bottom expression is given by $y = -h(x)$, where $h$ is the variable depth of the undisturbed water. The variable $y=\eta(x,t)$ represents the surface displacement, and the velocity components are $u(x,y,t)$ and $v(x,y,t)$.

{Incompressible Euler equations in a vertical $x$-$y$ plane, in the domain occupied by the fluid (see Figure \ref{shallow_water1})
 are given by }
\begin{eqnarray}
    u_t + u u_x + v u_y & = & -p_x/\rho \label{Euleru}\\
    v_t + u v_x + v v_y & = & -p_y/\rho - g \label{Eulerv}
\end{eqnarray}
The continuity equation is 
\begin{equation}
\label{eq_continuity}
    u_x + v_y = 0.
\end{equation}
{It can be shown that, because of the different scaling in $x$ and $y$, the vertical component $v$ of the velocity is much smaller than the horizontal one $u$. Furthermore, in incompressible inviscid flow, vorticity $\omega = v_y-u_x$ is carried unchanged aling material lines. As a consequence, if the vorticity is zero at initial time, it will remain zero for all time. As a consequence of these to facts, one can assume that $u_y$ is negligible, therefore in the shallow water approximation the horizontal velocity is independent of $y$. }
As a consequence, Eq.~\eqref{Euleru} becomes 
\[
    u_t + u u_x   =  g\eta_x(x,t)
\]
with $u = u(x,t)$.

At the free surface, the kinematic condition to be satisfied is
\begin{equation}
\label{eq_kincond}
    (\eta_t + u\eta_x - v) \Big|_{y=\eta}=0
\end{equation}
and the dynamic condition is 
\begin{equation}
\label{eq_dyncond}
    p \Big|_{y=\eta}=0.
\end{equation}

At the bottom, the only condition to be satisfied is 
\begin{equation}
\label{eq_botcond}
    (uh_x+v) \Big|_{y=-h}=0.
\end{equation}
Integrating Eq.~\eqref{eq_continuity} with respect to $y$, it becomes 
\begin{equation} 
\label{eq_integcontinuity}
    \int_{-h}^\eta(u_x)dy+v\Big|_{-h}^\eta=0.
\end{equation}
Substituting the conditions \ref{eq_kincond} and \ref{eq_botcond} in Eq.~\eqref{eq_integcontinuity}, it follows that
\begin{equation} 
\label{eq_relation}
    \int_{-h}^\eta u_xdy+\eta_t+u|_\eta\cdot\eta_x+u\Big|_{-h}\cdot h_x=0.
\end{equation}
Now, this relation is valid
\begin{equation}
   \frac{\partial}{\partial x} \int_{-h(x)}^{\eta(x)} udy = u\Big|_{y=\eta}\cdot\eta_x+u\Big|_{y=-h}\cdot h_x + \int_{-h}^\eta u_xdy
\end{equation}
and combining it with Eq.~\eqref{eq_integcontinuity}, it becomes
\begin{equation}
\label{eq_deriv}
    \frac{\partial}{\partial x} \int_{-h}^{\eta} udy = -\eta_t.
\end{equation}
Until here, the formulation is exact because no approximation has been applied. 

The shallow water approximation results from the assumption that the pressure $p$ is given by 
\begin{equation}
    p=g\rho(\eta-y)
\end{equation}
as in hydrostatics, where $\rho$ is the water density. In this way $p$ is assumed to be independent from the $y$-component of the acceleration of the water particles. 

Deriving the pressure with respect to $x$ gives
\begin{equation}
\label{eq_px}
    p_x=g\rho\eta_x
\end{equation}
and thus $p_x$ is independent of $y$. Therefore, it follows that the acceleration of the water particles in the $x$ direction is also independent of $y$, and so is $u$, the $x$-component of the velocity. So, from Eq. \ref{eq_px} the equation of motion in the $x$-direction becomes
\begin{equation}
\label{diff_eq_1}
    u_t+uu_x=-g\eta_x.
\end{equation}
This equation makes use of $u_y=0$, and Eq. \ref{eq_deriv} can be written as
\begin{equation}
\label{diff_eq_2}
    [u(\eta+h)]_x=-\eta_t
\end{equation}
having $\displaystyle \int^{\eta}_{-h} u dy = u \int^{\eta}_{-h} dy$.
Equations \ref{diff_eq_1} and \ref{diff_eq_2} are the first order differential equations that represent the non-linear shallow water theory. If the initial conditions for $u$ and $\eta$ are specified, the solution of these two equations describes the motion of the fluid.

\subsection{Laplace-Beltrami operator}
\label{A:LB}
{ 

Here we prove that 
\begin{equation}
    \label{A:LapBel}
    \partial_i\tpar_i c = \tpar_i\tpar_i c =: \Delta_\perp c
\end{equation}
where $\Delta_\perp$ is the Laplace-Beltrami operator. We also provide a representation of $\Delta_\perp$ in Cartesian coordinates.

First, observe that
\[
    \partial_i = \tpar_i  + n_in_j\partial_j,
\]
therefore 
\[
    \partial_i\tpar_i c = \tpar_i\tpar_i c + n_in_j\partial_j\tpar_i c.
\]
Now consider the term
\[
    n_in_j\partial_j\tpar_i c = n_i\pad{}{n}\tpar_i c.
\]
Since $n_i \tpar_i c \equiv 0$, it is
\[
    0 = \pad{}{n}(n_i\tpar_i c) = \pad{n_i}{n}\tpar_i c + n_i \pad{}{n}\tpar_i c.
\]
Now $\partial n_i/\partial n = n_j\partial_j n_i = 0$, therefore 
\[
    n_i\pad{}{n}\tpar_i c = 0
\]
from which \ref{A:LapBel} follows.

Now we give a representation for $\Delta_\perp$.    
We start observing that $\tpar_i c = h_{ij}\partial_j c$ is orthogonal to $\Sigma_B$, since
$n_ih_{ij} = 0$. 
It is
\begin{eqnarray}
    \Delta_\perp c = \tpar_i\tpar_i c &  =& \tpar_i(h_{ij}\partial_j c) = \tpar_i(\partial_i c - n_i n_j \partial_j c) \\
    & =& h_{ij}\partial^2_{ij} c - \tpar_in_i\pad{c}{n} - n_i (\tpar_i n_j) \partial_j c - n_i n_j \tpar_i\partial_j c
\end{eqnarray}
Observe that 
\[
    \tpar_i n_i = h_{ij}\chi_{ji} = \nabla\cdot \hat{n}
\]
where $\chi_{ij} \equiv \partial_i n_j$ is the {\em second fundamental form\/} of the surface in Cartesian coordinates $\Sigma_B$ and is orthogonal to $\hat{n}$. 
The last two terms are zero, because 
$n_i\chi_{ij} = 0$ and $n_i\tpar_i\partial_j c = 0$.

Therefore it follows:
\[
    \Delta_\perp c = \tpar_i\tpar_i c  = 
    h_{ij}\partial^2_{ij}c - \chi_{ii} n_j\partial_j c
\]
We recall that the trace of the second fundamental form of a 2D surface in 3D is the sum of the two principal curvatures, i.e. twice the mean curvature of the surface. 

\subsubsection{Laplace-Beltrami operator in two dimensions}
In two space dimensions (see Fig.~\ref{fig:bubble}) expressions considerably simplify.
Here we prove Eq.~\eqref{eq:BL2D}, i.e. 
\begin{equation}
    \frac{\partial^2 c}{\partial \tau^2} \equiv
    \tau_j\partial_j(\tau_i\partial_i c) = 
    h_{ij}\partial^2_{ij}c - \chi_{ii}\pad{c}{n} \equiv \Delta_\perp c.
    \label{A2:LB}
\end{equation}
Indeed it is 
\begin{equation}
    \label{A:tau2}
    \frac{\partial^2 c}{\partial \tau^2} = 
    \tau_j\partial_j(\tau_i\partial_i c) = 
    \tau_i\tau_j\partial^2_{ij} c + \tau_j(\partial_j\tau_i)\partial_i c = 
    \tau_i\tau_j\partial^2_{ij} c - \frac{1}{R} n_i\partial_i c
\end{equation}
the latter equality is a consequence of the relation 
\[
    \pad{\hat{\tau}}{\tau} = -\frac{1}{R} \hat{n}
\]
where $R$ denotes the local radius of curvature.
Finally we observe that in 2D it is 
\begin{equation}
    \label{A:h}
    h_{ij} \equiv \delta_{ij}-n_in_j = \tau_i\tau_j
\end{equation}
as it can be easily shown by inspection checking that 
\[
    n_i n_j + \tau_i \tau_j = \delta_{ij}
\]
Using Eq.~\eqref{A:tau2} and \ref{A:h} in \ref{A2:LB} completes the proof of relation \ref{eq:BL2D}.
}

\begin{exercises}
\item Implement the discretization of the Eqs.~(\ref{reduced1d}-\ref{BCt}) in $\Omega^0 = [0,1]$ using second order accurate Finite Differences scheme (see Chapter ... for more details). Note that Eq.~\eqref{BCt} is already discretized in Eq.~\eqref{eq_BCt_discr}.
\begin{figure}[tb]
	\begin{minipage}
		{.44\textwidth}
		\centering
\begin{overpic}[abs,width=\textwidth,unit=1mm,scale=.25]{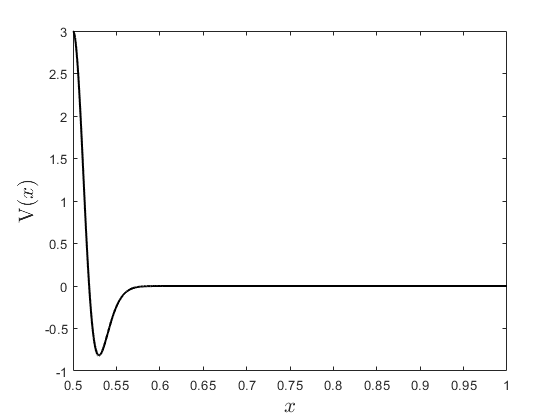}
\put(1,40){(a)}
\end{overpic}	\end{minipage}
	\begin{minipage}
		{.44\textwidth}
\begin{overpic}[abs,width=\textwidth,unit=1mm,scale=.25]{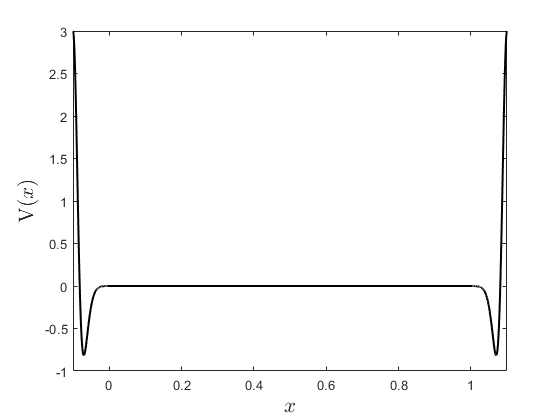}
\put(1,40){(b)}	\end{overpic}
\end{minipage}
	\caption{\textit{(a): Plot of the functional defined in Eq.~\eqref{eq_functional_exp}. (b): Plot of the functional defined in Eqs.(\ref{eq_functional_exp_2}-\ref{eq_functional_exp_f}).}}
 \label{fig_V_exerc}
\end{figure}
\item Now consider a 2D domain, in which the potential is defined as a function $W:=W(x,y)$ and it is constant in $y$-direction, i.e., $W(x,y) = V(x)$ (see Eq.~\eqref{expr_V_LJ}). In this case the domain is defined as $\Omega^\varepsilon = [-\varepsilon,1]\times[0,1]$. Discretize in space the following problem with Finite Differences scheme
\begin{align}
	\frac{\partial c}{\partial t} & = D\nabla\cdot \vec J \quad {\rm in }\, \Omega^\varepsilon\\
 \vec J &=  \nabla c + \frac{1}{k_BT}  c  \nabla W \\
\vec J \cdot \vec n & = 0 \quad {\rm on} \, \partial \Omega^\varepsilon 
\end{align}
where $\vec n$ is the outgoing normal vector to $\Omega^\varepsilon$.
\item Starting from the system~(\ref{reduced1d}-\ref{BCt})
\begin{align}
	\displaystyle \frac{\partial   c }{\partial t} + \frac{\partial  J}{ \partial x} & = 0 \qquad {\rm in} \, \Omega^\varepsilon  \\
 	\displaystyle J & = - D\left( \frac{\partial  c  }{ \partial x} + \frac{1}{k_BT}  c  V'\right)\\
  J &= 0 \qquad {\rm on} \, \partial\Omega^\varepsilon 
\end{align}
where, in this case, the full domain is $\Omega = [0.5,1]$, with $\Omega^\veps_b = [0.5,0.5+L\veps]$ the potential domain, $\Omega^\varepsilon_f = [0.5+L\varepsilon,1]$ the fluid domain, and $L = 2$.  

Following the strategy in Section~\ref{section_1Dmodel}, find the new expression for $M$ (see Eq.~\eqref{expr_M}) considering the following expression for the potential (see Fig.~\ref{fig_V_exerc} (a))
\begin{equation}
\label{eq_functional_exp}
  V(x) = a_1\exp\left(-b_1\left(\frac{x - x_0}{\varepsilon}\right)^2 \right) - a_2\exp\left(-b_2\left(\frac{x - x_0}{\varepsilon}\right)^2 \right)  
\end{equation} 
with $a_1 = 6, b_1 = 30, a_2 = 3, b_2 = 10, x_0 = 0.5$.

\item Starting from the system~(\ref{reduced1d}-\ref{BCt})
\begin{align}
	\displaystyle \frac{\partial   c }{\partial t} + \frac{\partial  J}{ \partial x} & = 0 \qquad {\rm in} \, \Omega^\varepsilon  \\
 	\displaystyle J & = - D\left( \frac{\partial  c  }{ \partial x} + \frac{1}{k_BT}  c  V'\right)  \\
  J &= 0 \qquad {\rm on} \, \partial\Omega^\varepsilon 
\end{align}
where, in this case, the potential $V$ is defined as follows (see Fig.~\ref{fig_V_exerc} (b))
\begin{align}
\label{eq_functional_exp_2}
V(x) = V_1(x) + V_2(x) \\
        V_1(x) = f(x,x_1), \quad V_2(x) = f(x,x_2), \quad x_1 = 0, \quad x_2 = 1,
\end{align}
and 
\begin{equation}
\label{eq_functional_exp_f}
f(x,x_l) = a_1\exp\left(-b_1\left(\frac{x - x_l}{\varepsilon}\right)^2 \right) - a_2\exp\left(-b_2\left(\frac{x - x_l}{\varepsilon}\right)^2 \right), 
\end{equation}
with $l = 1,2$, $a_1 = 6, b_1 = 30, a_2 = 3, b_2 = 10, x_0 = 0.5$.
The full domain is $\Omega = [-\varepsilon,1+\varepsilon]$, with $\Omega^\veps_{1,b} = \varepsilon[-1,L]$ and $\Omega^\veps_{2,b} = [1-L\varepsilon,1+\varepsilon]$ the potential domains, $\Omega^\varepsilon_f = [L\varepsilon,1-L\varepsilon]$ the fluid domain, and $L = 2$.  

Following the strategy in Section~\ref{section_1Dmodel}, find the new boundary conditions (see Eqs.~(\ref{BCt}-\ref{expr_M})).

\begin{figure}
\centering
\begin{minipage}
		{.44\textwidth}
		\centering
\begin{overpic}[abs,width=\textwidth,unit=1mm,scale=.25]{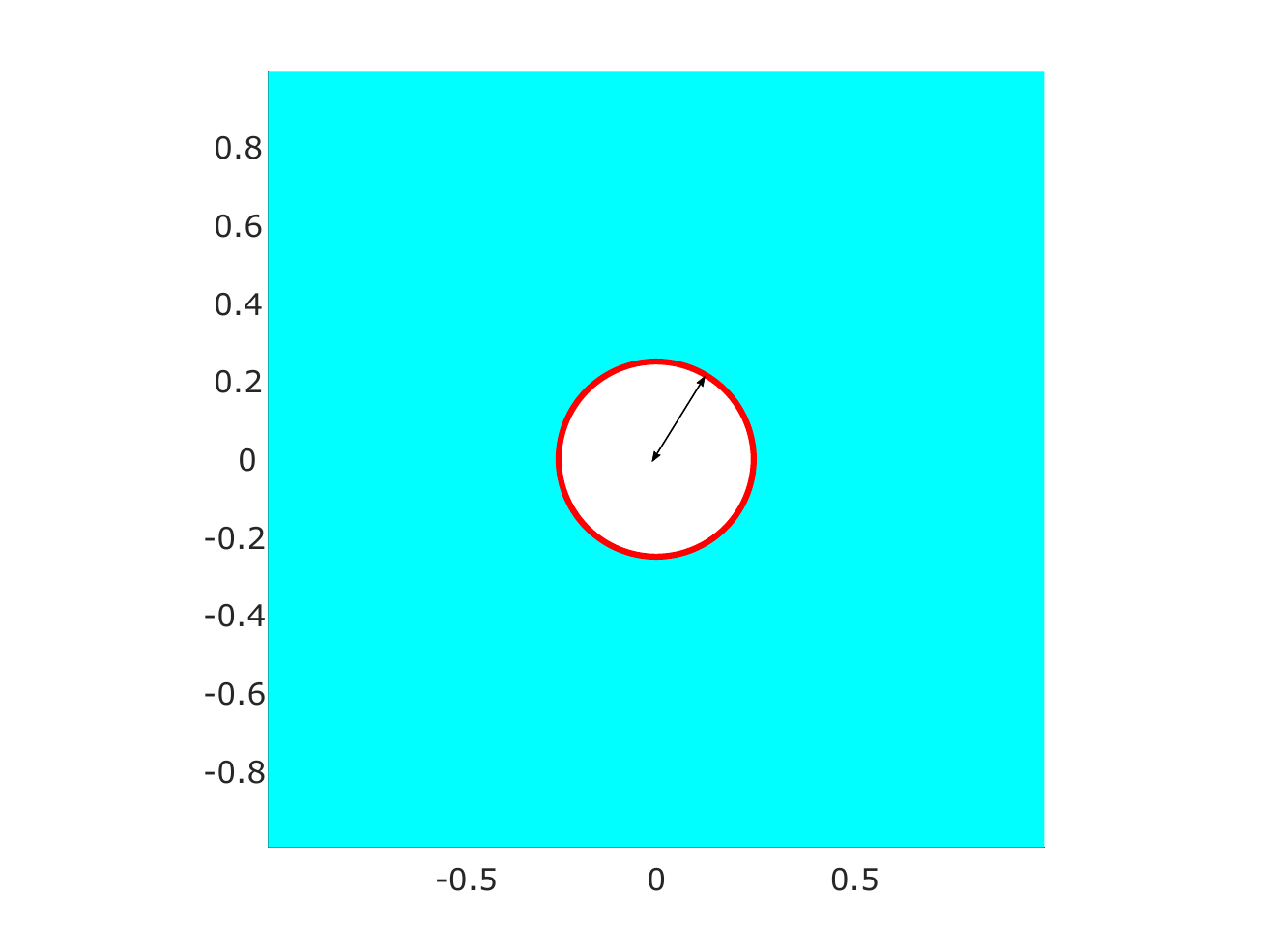}
\put(1,40){(a)}
\put(15,30){$\Omega^0$}
\put(30,33){$ \varphi < 0$}
\put(29.5,17){$\Big\uparrow$}
\put(27,12){$\varphi > 0$}
\put(27,26){$R_b$}
\end{overpic}
\end{minipage}
	\begin{minipage}
		{.44\textwidth}
\begin{overpic}[abs,width=\textwidth,unit=1mm,scale=.25]{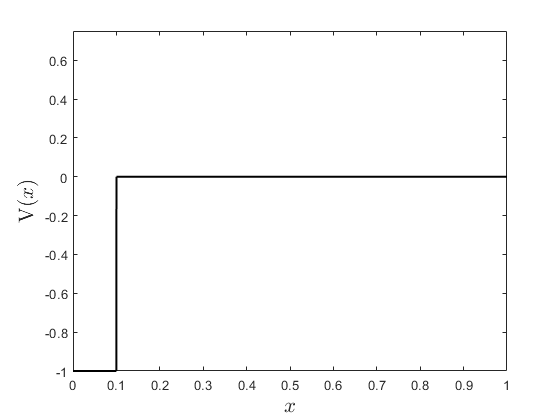}
\put(1,40){(b)}
\end{overpic}
\end{minipage}
\caption{\textit{(a): Scheme of the level-set function. (b): Scheme of a potential V(x) that represents a squared sink.}}
	\label{fig_domain_level_set}
\end{figure}

\item Let us consider the \textit{multiscale model} in Eqs.(\ref{system_multiscale}-\ref{bc_multiscale}), with $\Omega^0 = [-1,1]^2\setminus \mathcal{B}$ and a ball $\mathcal{B}$ centered in $C=(0,0)$, with radius $R_b = 0.25$. Using a level-set approach, it is possible to implicitly know the  boundaries of the ball (that are marked in red in Fig.~\ref{fig_domain_level_set} (a)) defining a function $\varphi = R_{b}-\sqrt{x^2+y^2}, \, (x,y) \in [-1,1]^2$. With this technique, the values of the level-set function are negative inside the domain and positive outside of the domain (i.e., the points that are inside the ball; see Fig.~\ref{fig_domain_level_set} (a)). 

Given two balls in the domain $\Omega^0$, s.t. $\Omega^0 = [-1,1]^2\setminus(\mathcal{B}_1\cup\mathcal{B}_2)$, with $\mathcal{B}_1$ centered in $C_1 = (-0.5,-0.5)$ and radius $R_{b,1} = 0.2$, and $\mathcal{B}_2$ centered in $C_2 = (0.5,0.5)$ and radius $R_{b,2} = 0.25$. Determine a possible level set function $\varphi$, s.t. $\varphi > 0 $ in $\mathcal{B}_1 \cup \mathcal{B}_2$, and $\varphi < 0 $ in $\Omega^0$.

\item Given the initial conditions in Eq.~\eqref{ic_trigonom}, determine the constants $a,b$ in Eq.~\eqref{eq_general_solu}.
\item  Given the initial conditions in Eq.~\eqref{ic_trigonom}, determine the constants $C_+,C_-$ in Eq.~\eqref{eq_gener_solu_complex}.

\item Find the expression of the solution of Eq.~\eqref{eq_lambda_2nd} in the three case $\Delta >0,\Delta =0, \Delta <0$.

\item \label{ex:pendulum} Deduce the equation of motion of the pendulum starting from equations \eqref{eq:pendulum2} and \eqref{eq:x1_x2}.

\item Prove that Eq.~\eqref{eq_TFT} is valid.
\end{exercises}

\begin{project}
Find an analytical solution in 1D of the \textit{multiscale model} in Eqs.~(\ref{reduced1d}-\ref{BCt}), with a squared sink in $x = 0$, as it shows in Fig.~\ref{fig_domain_level_set} (b).
\end{project}
\begin{project}
In a 2D squared domain $\Omega^0 = [-1,1]^2\setminus \mathcal{B}$, solve numerically the system~(\ref{system_multiscale}-\ref{bc_multiscale}) with a squared hole $\mathcal{B}$ of coordinates: $A=(-0.5,-0.5), \, B=(-0.5,0.5), \, C = (0.5,0.5), \, D = (0.5,-0.5)$.   
\end{project}

\backmatter

\cleardoublepage

\printindex

\end{document}